\documentclass[twoside,12pt,draft]{article}

\input{amssym.def}
\input{amssym.tex}

\newtheorem{theorem}{Theorem}[section]
\newtheorem{proposition}[theorem]{Proposition}
\newtheorem{lemma}[theorem]{Lemma}
\newtheorem{definition}[theorem]{Definition}
\newtheorem{example}[theorem]{Example}
\newtheorem{remark}[theorem]{Remark}

\newcommand{\btheorem}{\begin{theorem}}
\newcommand{\etheorem}{\end{theorem}}
\newcommand{\bprop}{\begin{proposition}}
\newcommand{\eprop}{\end{proposition}}
\newcommand{\blemma}{\begin{lemma}}
\newcommand{\elemma}{\end{lemma}}
\newcommand{\bdef}{\begin{definition}\rm}
\newcommand{\bex}{\begin{example}\rm}
\newcommand{\eex}{\end{example}}
\newcommand{\bremark}{\begin{remark}\rm}
\newcommand{\eremark}{\end{remark}}

\newcommand{\qed}{\ \hfill \rule{2mm}{2mm}}

\newenvironment{proof}{\noindent PROOF:~}{\qed\medskip\par} 
\newcommand{\bproof}{\begin{proof}}
\newcommand{\eproof}{\end{proof} \hfill \vspace{6mm}}

\newcommand{\bvgln}{$$\begin{array}{lcl}}
\newcommand{\evgln}{\end{array}$$}
\newcommand{\bvgl}{\begin{eqnarray*}}
\newcommand{\evgl}{\end{eqnarray*}}

\newcommand{\comult}{\Phi}
\newcommand{\haar}{h}
\newcommand{\antipode}{\kappa}
\newcommand{\cst}{C\/$^{\ast}$}
\newcommand{\st}{$^{\ast}$}
\newcommand{\ten}{\otimes}

\newcommand{\R}{\mbox{$\Bbb R$} }
\newcommand{\C}{\mbox{$\Bbb C$} }
\newcommand{\eps}{\varepsilon}
\renewcommand{\to}{\rightarrow}
\newcommand{\tr}{\mbox{Tr}}
\newcommand{\tenrep}{\mbox{ $\mbox{\scriptsize \sf T}
\hspace{-1.8ex}\bigcirc$}}

\renewcommand{\H}{\mbox{$\mathcal{H}$}}
\newcommand{\K}{\mbox{$\mathcal{K}$}}
\newcommand{\B}{\mbox{$\mathcal{B}$}}

\begin{document}

\title{Notes on Compact Quantum Groups}

\author{Ann Maes\thanks{Research Assistent of the
National Fund for Scientific Research (Belgium)}
\space \& Alfons Van Daele
\\Department of Mathematics
\\Katholieke Universiteit Leuven \\Celestijnenlaan 200B \\B-3001 Heverlee
\\Belgium \\email: Ann.Maes@wis.kuleuven.ac.be  \\
Alfons.VanDaele@wis.kuleuven.ac.be}

\date{Maart 1998}

\maketitle

\small
\subsection*{Abstract}
Compact quantum groups have been studied by several authors and from
different points of view. The difference lies mainly in the choice of
the axioms. In the end, the main results turn out to be the same.
Nevertheless, the starting point has a strong influence on how the main
results are obtained and on showing that certain examples satisfy these
axioms.
\smallskip

In these notes, we follow the approach of Woronowicz and we treat the
compact quantum groups in the C$^*$-algebra framework. This is a natural
choice when compact quantum groups are seen as a special case of locally
compact quantum groups. A deep understanding of compact quantum groups
in this setting is very important for understanding the problems that
arise in developing a theory for locally compact quantum groups.
\smallskip

We start with a discussion on locally compact quantum
groups in general but mainly to motivate the choice of the axioms for the
compact quantum groups. Then we develop the theory. We give the main
examples and we show how they fit into this framework.
\smallskip

The paper is an expository paper. It does not contain many new results
although some of the proofs are certainly new and more elegant than  the
existing ones. Moreover, we have chosen to give a rather complete and
self-contained treatment so that the paper can also serve as an
introductory paper for non-specialists. Different aspects can be learned
from these notes and a great deal of insight can be obtained.

\normalsize

\section{Introduction}

Quantum groups have been widely studied the last decade.  Most of all,
they were investigated in a purely algebraic context and examples played
an important role.  Not so much has been done for a comprehensive theory
of locally compact quantum groups.  By now, the reason for this has
become clear.  It is twofold.  First of all, it seems to be difficult,
both technically and conceptually, to develop such a theory.
Secondly, examples that are simple on the algebraic level turn out to get
sometimes very complicated when lifted to topological level. Nevertheless
there are many good reasons to develop such a theory of locally compact
quantum groups.  The main advantage of such a theory is the
possibility of using the equivalent of the Haar measure for locally
compact groups.

In section 2 of these notes, we say a little more about the ideas behind
a possible theory of locally compact quantum groups.  But since this
is not the main purpose of these notes, we do not enter very deeply
into it.   General theories of locally compact quantum groups exist
(see e.g. \cite{baaj, es, m, mn, mnw}), but it must be said that they
are not yet completely satisfactory.  A lot of research still 
needs to be done there.

The situation is quite different for discrete and compact quantum groups.
These two classes are very well understood.  In fact there is now
also a class of locally compact quantum groups, that contains the discrete
and the compact ones, and is self dual in the sense that the
Pontryagin dual of a locally compact quantum group of this type is again
of this type.  These quantum groups are studied in \cite{vd.alg, vd.alg2,
kustvd}.

The compact quantum groups are certainly those locally compact quantum 
groups that have been studied most.  We just mention some of the 
references here: Woronowicz~\cite{wor1, wor2, wor3}, Kirchberg~\cite{kir},
Dijkhuizen and Koornwinder~\cite{dijk}.  These are the places where compact
quantum groups were developed.  But there are also other authors that
treated compact quantum groups in duality with discrete quantum groups.
This is e.g. the case in the paper~\cite{baaj} by Baaj and Skandalis.  
And remark that the paper on compact
quantum groups \cite{dijk} by Dijkhuizen and Koornwinder is very close to
the paper \cite{eff} by Effros and Ruan on discrete quantum groups.  There has indeed been 
some ambiguity.  For instance, Kirchberg called his objects discrete
quantum groups while others would speak about compact quantum groups.

When looking at compact quantum groups as compact, locally compact 
quantum groups, the obvious approach is the one of Woronowicz.  He first
developed the theory of compact matrix quantum groups in~\cite{wor1}.
Later, an important note \cite{wor2} was published where the axioms could be 
simplified.  In 1992, a preprint, entitled 'Compact Quantum Groups'
\cite{wor3} was widely distributed, but it was not completed for a long time.
Now, it seems that the paper is going to appear.

Because of this ambiguity and the different approaches by different authors, 
it seems desirable to have a text which starts from a natural set of axioms 
and developes the theory in a fairly complete way.  This is what we do 
in this paper.  We work with the definition of Woronowicz since we view 
compact quantum groups as compact locally compact quantum groups, and 
therefore, the \cst-approach is the natural choice.  We start by motivating
this choice.  Then we develop the theory as it can be found in different
articles.  To a great extend however, we follow the general treatment
of Woronowicz in~\cite{wor3}.  But we also explain the link with the
earlier paper on compact matrix pseudogroups~\cite{wor1} and the 
note~\cite{wor2}, simplifying the axioms for such compact matrix 
pseudogroups.  Moreover, we include the important examples in the 
beginning.  And we end the paper by treating the Pontryagin dual 
of a compact quantum group using more recent ideas.

This paper is expository.  The results are not new.  Nevertheless, some of 
the proofs, also of the more important ones, are original.  Some novelty
also lies in the structure of our treatment.  In any case, we have made
an effort to write down notes that are self-contained to a great extend.  
This should help the reader who wants to study this approach to compact 
quantum groups without the need to find his way in the different articles on
the subject (with different approaches).

The paper is organized as follows.  In section~2 we introduce the concept of
quantization, as a way to generalize locally compact groups.  We describe how
this process works in the case of finite quantum groups, and
what the problems are when we want to extend this construction to the locally
compact case.  We then motivate how Woronowicz' treatment of compact quantum
groups fits into this scheme.
In section~3 we define a compact quantum group.  We motivate and discuss this
definition and we give some examples.  We also introduce the compact matrix
quantum groups defined in~\cite{wor2}
and show that they indeed satisfy this set of axioms.
One of the main features of this definition is that, starting from a simple set of
axioms, very strong results can be proven. Such a result is the existence and
the uniqueness of the Haar measure.  This is proven and discussed in section~4.

In section~5, the concept of representation is translated to the quantum
language.  Then, like in the classical group case, the Haar measure is used to
construct the right regular representation.
In section~6, we go deeper into the representation theory.  The main results in
this section are that any unitary representation can be decomposed into
irreducible unitary representations and that all irreducible unitary
representations are contained in the regular representation.

In section~7, we look at the subspace spanned by the matrix elements of 
finite-dimensional unitary representations.  We show that it is a dense
\st-subalgebra and that it can be made into a Hopf \st-algebra.
It then becomes very straightforward to prove that Woronowicz' new definition
of a compact matrix quantum group in\cite{wor2} is equivalent with the
original one in~\cite{wor1}.
In the last section we generalize the Pontryagin duality between compact 
and discrete abelian groups to a duality theory between compact quantum 
groups and discrete quantum groups.

For the theory of \cst-algebras, we refer to the book of Dixmier~\cite{dix}.
For general results
on Hopf algebras, we refer to the books of Abe~\cite{abe} and 
Sweedler~\cite{swe}, and for Hopf
$ ^{\ast}$-algebras and the duality between them, to \cite{vd.dual}.
For any Hilbert space $\H$, the \cst-algebra of bounded linear operators 
on $\H$ will be denoted by $\B(\H)$ and the \cst-algebra of compact 
operators by $\B_0(\H)$.  For given vectors $\xi,\eta \in \H$, we will 
denote by $\omega_{\xi,\eta}$ the linear functional on $\B(\H)$ such 
that $\omega_{\xi,\eta}(x)=\langle x\xi,\eta \rangle$ for all $x \in \B(\H)$.  
If $A$ and $B$ are \cst-algebras, then the algebraic tensor product will
be denoted by $A \odot B$.  We will always use the minimal \cst-norm and
the completion of $A \odot B$ with respect to this norm will be denoted
by $A \ten B$.  We will use the notation $\iota$ for the identity map.
If $\omega$ is a continuous linear functional on $A$, one 
can consider the linear mapping $(\omega \ten \iota):A \odot B \to B$ 
defined by $(\omega \ten \iota)(a \ten b) = \omega(a)b$ if $a \in A$ and
$b \in B$.  This mapping is continuous with respect to the minimal \cst-norm
 on $A \odot B$, and hence can be extended to a linear mapping $A \ten B
\to B$, which will still be denoted by $\omega \ten \iota$.  In the same way 
$\iota \ten \omega$ is defined.  These maps are called slice maps.

\section{Locally compact quantum groups}

Let $A$ be a \cst-algebra.  Let $A \ten A$ be the \cst-tensor product of $A$
with itself with respect to the minimal \cst-tensor product norm.  Let
$M(A \ten A)$ be the multiplier algebra of $A \ten A$.   A \st-homomorphism
$\comult : A \to M(A \ten A)$ is called non-degenerate if $\comult (A)(A \ten
A)$ is dense in $A \ten A$.  Denote by $\iota$ the identity map on $A$ and
consider the \st-homomorphisms $\comult \ten \iota$ and $\iota \ten \comult$
on the algebraic tensor product $A \odot A$ of $A$ with itself.  By
continuity, they have unique extensions to \st-homomorphisms from $A \ten A$
to $M(A \ten A \ten A)$.  If $\comult$ is non-degenerate, there are also
unique extensions to \st-homomorphisms from $M(A \ten A)$ to $M(A \ten A
\ten A)$ (see e.g. \cite{tom}).
These extensions will still be denoted by $\comult \ten \iota$
and $\iota \ten \comult$.

\bdef
A non-degenerate \st-homomorphism $\comult : A \to M(A \ten A)$ is called
coassociative if $(\comult \ten \iota) \comult = (\iota \ten \comult) \comult$.
Then we call $\comult$ a comultiplication on $A$.
\end{definition}

Let $G$ be a locally compact space.  Let $A$ be the \cst-algebra $C_0(G)$
of continuous complex functions on $G$, vanishing at infinity.  Here, $A
\ten A$ can be identified with $C_0(G \times G)$ and $M(A \ten A)$ with
the \cst-algebra $C_b(G \times G)$ of bounded complex functions on
$G \times G$.  If $G$ is endowed with a continuous, associative
multiplication, then we can define a comultiplication $\comult$ on $A$
by $(\comult (f))(p,q) = f(pq)$ whenever $f \in C_0(G)$ and $p,q \in G$.  It
is not hard to see that $\comult$ is indeed non-degenerate.  And of course,
the coassociativity follows from the associativity of the product on $G$.
So, with any locally compact semi-group, we can associate, in a canonical
way, a pair $(A, \comult)$ of an abelian \cst-algebra $A$ and a
comultiplication $\comult$ on $A$.

Conversely, suppose that $A$ is an
abelian \cst-algebra and $\comult$ a comultiplication on $A$, then $A$ has
the form $C_0(G)$ for some locally compact space $G$ and $\comult$ defines
a multiplication on $G$ by the formula $f(pq) = (\comult (f))(p,q)$.  The
non-degeneracy of $\comult$ must be used to do this.  And again, the
multiplication is associative because $\comult$ is coassociative.  Hence $G$
is a locally compact semi-group.

This shows that there is a one-to-one correspondence between the locally
compact semi-groups on the one hand and the pairs $(A, \comult)$ of an
abelian \cst-algebra $A$ with a comultiplication $\comult$ on the other hand.
The \cst-algebra structure takes care of the underlying topology, while the
comultiplication reflects the multiplication.

This brings us to the following notion.

\bdef
A locally compact quantum semi-group is a pair $(A, \comult)$ of a
\cst-algebra $A$ and a comultiplication $\comult$ on $A$.
\end{definition}

The next natural question is~: What is a locally compact quantum group~?
Since a locally compact group is a locally compact semi-group, where the
multiplication satisfies some extra conditions, it is most natural that
a locally compact quantum group will be a pair $(A, \comult)$ of a
\cst-algebra $A$ with a comultiplication $\comult$ on $A$ satisfying some
extra properties.  What are such conditions~?

To distinguish a group from a semi-group, the existence of a unit and of
an inverse for all the elements is required.  If $G$
is a locally compact group, the identity gives rise to a \st-homomorphism
$\eps : C_0(G) \to \C$ given by $\eps (f) = f(e)$.  The inverse gives rise
to an anti-homomorphism $\antipode : C_0(G) \rightarrow C_0(G)$ given by
$(\antipode(f))(p) =
f(p^{-1})$ whenever $p \in G$.  The properties of the identity and the
inverse can easily be translated into properties of $\eps$ and $\antipode$.
However, as it turns out, to obtain the quantum analogues of the maps
$\eps$ and $\antipode$ and their natural properties seems to be rather complicated.

We will explain these difficulties by looking first at the case of a finite
group where this whole procedure presents no problems at all.

If $G$ is a finite group, the \cst-algebra $C_0(G)$ is the \st-algebra
$C(G)$ of all complex functions on $G$ with pointwise operations.
It is of course finite-dimensional
so that it suffices to look at the algebraic tensor product $C(G) \odot
C(G)$. Also the multiplier algebra is this algebraic tensor product.
In particular, the comultiplication $\comult$ is a comultiplication in the
ordinary, algebraic sense on the algebra $C(G)$.

The properties $ep = pe = p$ and $p^{-1}p = pp^{-1} =e$ for all $p$, become
in terms of $\eps$ and $\antipode$
$$(\eps \ten \iota)\comult(f) = (\iota \ten \eps)\comult(f) = f$$
$$m(\antipode \ten \iota) \comult(f) = m(\iota \ten \antipode)\comult(f) = f$$
where $\iota$ denotes the identity map as before, and where $m$ is the
multiplication as a map from the tensor product $C(G) \odot C(G)$ to $C(G)$
defined by $m(f \ten g) = fg$.
This means that $(C(G), \comult)$ is a Hopf algebra.  In fact, it is a Hopf
\st-algebra if we let $f^{\ast}(p) = \overline{f(p)}$.

Let us recall the definition of a Hopf \st-algebra
(see e.g. \cite{vd.dual})~:

\bdef
Let $A$ be an associative \st-algebra over \C with an identity.  Assume
that $\comult$ is a \st-homomorphism from $A$ to the algebraic tensor product
$A \odot A$ satisfying $\comult (1) = 1 \ten 1$ and coassociativity $(\comult
\ten \iota) \comult = (\iota \ten \comult) \comult$.  Also assume that there are
linear maps $\eps : A \to \C$ and $\antipode : A \to A$ such that
$$(\eps \ten \iota)\comult(a) = (\iota \ten \eps) \comult(a) = a$$
$$m(\antipode \ten \iota) \comult(a) = m(\iota \ten \antipode)\comult(a) = \eps(a)1$$
for all $a \in A$.
Then $(A, \comult)$ is called a Hopf \st-algebra.
\end{definition}

If they exist as linear maps satisfying the above equations,
the maps $\eps$ and $\antipode$ are unique.  We have automatically
that $\eps$ is a \st-homomorphism, that $\antipode$ is an anti-homomorphism and that
$\antipode(\antipode(a)^{\ast})^{\ast} = a$ (see e.g.~\cite{vd.dual}).

So, in this setting, it is natural to define a finite quantum group in
the following way~:

\bdef \label{def2.4}
A finite quantum group is a pair $(A, \comult)$ of a finite-dimensional
\cst-algebra $A$ with a comultiplication $\comult$ such that $(A, \comult)$
is a Hopf \st-algebra.
\end{definition}

Of course, there are other possible choices.  The most general notion would
be to take any Hopf algebra $(A, \comult)$ with $A$ finite-dimensional.  If
we restrict to the case of Hopf \st-algebras with a finite-dimensional
\cst-algebra for the underlying algebra, we have the property that a
finite quantum group is a group if and only if the underlying algebra $A$ is
abelian.  We will see later in these notes that there are more reasons to
define finite quantum groups in this way (within the theory of locally
compact quantum groups).

Now we can discuss the main problems that arise if we want to generalize
this.  First of all, it turns out that the natural candidates for the maps
$\eps$ and $\antipode$ in interesting examples, are only densely defined, and not
continuous.  The counit $\eps$ is a \st-homomorphism on a dense
\st-subalgebra.  The fact that $\antipode$ is an anti-homomorphism, related with the
involution by the formula $\antipode(\antipode(a)^{\ast})^{\ast} = a$ causes some extra
difficulties.  Then it is no longer obvious how one has to define the maps
$\eps \ten \iota$, $\iota \ten \eps$, $\antipode \ten \iota$ and $\iota \ten \antipode$
on the image of $\comult$ which lies in $M(A \ten A)$.  A last difficulty
is the multiplication map $m$.  It is only defined on the algebraic tensor
product $A \odot A$ and it is very badly behaved with respect to both the
algebraic and the topological structure.

That these maps (the counit $\eps$, the antipode $\antipode$ and the multiplication
$m$) are badly behaved in the general quantum case, is certainly the main
reason why, at this moment, there is not yet a satisfactory theory of
locally compact quantum groups (in general).

Then, one might wonder why such a theory is desirable anyway.  What kind of
properties do we hope to obtain~?
We would like to have conditions on $\comult$ such that the
following program can be carried out.

First of all, if $(A, \comult)$ comes from a locally compact group as before,
then $\comult$ should satisfy the conditions.  On the other hand, if $A$ is
abelian and $\comult$ satisfies these conditions, one should be back in the
group case.  Remark that, in the finite case, with the
definition~\ref{def2.4}
for a finite quantum group, this is the case.  It should also be possible
to develop a nice theory for these locally compact quantum groups.  For
instance we would like to prove the existence and uniqueness of the Haar
measure, to have a nice representation theory...  We are also very much
interested in a duality theory extending the Pontryagin duality for
locally compact abelian groups.  It is expected that the dual of a locally
compact quantum group will again be such a locally compact quantum group.
Finally, the conditions on $\comult$ may not be too strong so that
sufficiently many examples fit into this frame.

As we mentioned already in the introduction, there exist notions of
locally compact quantum groups.  First, we have the Kac algebras 
(see~\cite{es}).  They have recently been generalized by Masuda and 
Nakagami~\cite{mn} so as to include also the more recent examples where
the antipode $\kappa$ is not a $\ast$-map.  These two theories are
formulated in the von Neumann algebra context.  It would be
desirable (and natural) to have a theory in the \cst-algebra framework.
The work of Masuda, Nakagami and Woronowicz~\cite{mnw} has been
announced at several 
conferences, among them in Warsaw~\cite{m}, but the paper is still not
available.

Generally, it is the opinion of many researchers in the field that 
there is still a long way to go.  But equally strong is the feeling that,
since there are several special cases and also nice examples not 
belonging to these special cases, it should be possible to come up
with a comprehensive theory in the future.


Now, what are the special cases that have been studied and are understood?
First there is of course the finite case.  It is known that for any
finite-dimensional Hopf algebra, there exist a left and a right
invariant integral
(the equivalent of the Haar measure).  With our definition of a finite
quantum group, the left and right invariant integrals coincide and are
positive (like what you expect in the finite case).
Of course, there is also a duality (see e.g.~\cite{vd.dual}).  The dual algebra is
again a finite-dimensional \cst-algebra.  This is an extra motivation for
our definition of a finite quantum group.  For a simple treatment of
finite quantum groups, we refer to~\cite{vd.fin}.

There are two other special cases that are well understood.  The first one
is the compact case and the second one is the discrete case.  Compactness
and discreteness are topological properties.  Therefore, they are reflected
into properties of the underlying \cst-algebra.

In the compact case, this is easy.  If $G$ is a compact group, then $C_0(G)$
coincides with the \cst-algebra of all complex continuous functions on $G$.
So $C_0(G)$ has an identity.  Therefore, a compact quantum group will be
a pair $(A,\comult)$ of a \cst-algebra $A$ with an identity and a
comultiplication $\comult$ satisfying certain extra conditions.  In section~3
of these notes we will give these conditions and argue about them.  And as
we will see further in these notes, this notion of a compact quantum group
yields a satisfactory theory with nice examples.

In the discrete case, the choice is not so obvious.  If $G$ is discrete,
then the multiplier algebra $C_b(G)$ of $C_0(G)$ coincides with the double
dual of $C_0(G)$.  The \cst-algebra of compact operators on a Hilbert space
also has this property, and direct sums of these.  There are however
several reasons to consider only direct sums of full matrix algebras
$M_n(\C)$ (the compact operators on a finite-dimensional Hilbert space).
So, a discrete quantum group is a pair $(A, \comult)$ of a \cst-algebra $A$
which is a \cst-direct sum of full matrix algebras satisfying certain
properties.  See~\cite{vd.discr} for a possible more precise definition.
Also here, there is a nice theory with interesting examples.

The compact case is dual to the discrete case in a sense that can be made
precise.  Using the Haar measure on a compact quantum group,
one can construct the dual, and this will be a
discrete quantum group.  We will show this in section~8 of these notes.
There we will also give some more information about the notion of a discrete
quantum group.  Conversely, starting with a discrete quantum group, one can
recover (more or less) the original compact quantum group by a similar
process.

Many of the different approaches to compact quantum groups are in one way
or another related with this duality between compact and discrete quantum
groups.  In \cite{kir} Kirchberg called a compact quantum group in our sense a
discrete quantum group.  And in \cite{eff} Effros and Ruan introduce discrete
quantum groups while other people would claim that in fact, they study
compact quantum groups from another point of view.  Similar is the approach
of \mbox{Dijkhuizen} and Koornwinder in \cite{dijk}. In most of these other
approaches, conditions are
formulated, not only on the algebra $A$ and the comultiplication $\comult$,
but also on the dual objects.  Very often then, it is easier to get to the
main results.  This is not at all surprising~: Part of the main properties
are then already included in the axioms.

We will discuss these different aproaches further in these notes.  Our
philosophy should be clear already.  A locally compact group is a locally
compact space with a group structure, compatible with the topological
structure.  So a locally compact quantum group is a locally compact
quantum space (a \cst-algebra) with a quantum group structure (a comultiplication
with certain properties).  Furthermore compactness is a property of the
underlying topological space and therefore reflected in a property of the
\cst-algebra (and not in extra properties on $\comult$).  And similarly
for discreteness.

In a discrete space, there is very little topology.  Therefore it is of
course not hard (and maybe more natural) to treat the discrete quantum
groups outside the \cst-framework.  And because of this ambiguity, due
to the duality between the compact and the discrete case, it is
also quite natural to treat the compact case without the use of
\cst-algebra theory, such as in \cite{dijk}.
A more general theory is set up in 
\cite{vd.alg,vd.alg2}, where duality is studied within the category
of Multiplier Hopf Algebras.  This framework includes both the discrete
and the compact quantum groups and the duality in this category includes
the duality between discrete and compact quantum groups.
This whole setting is purely algebraic, but can be fitted into the 
\cst-algebraic framework~:
In \cite{kust} and \cite{kustvd},
starting from such an algebraic quantum group,
a \cst-algebraic group in the sense of \cite{mnw} is constructed.

In these notes, we look at the compact and the discrete cases
as special cases of the more general locally compact quantum groups, and
hence it is more obvious to do this within the \cst-algebra framework.
In fact, one could say that the \cst-approach to quantum groups is the
study of locally compact quantum groups.

This approach to quantum groups is necessary if we want to look at the
Haar measure and if we want to study unitary representations for
non-compact quantum groups.  If e.g. we consider a non-compact matrix
group, the Haar measure will not be defined on the function algebra
generated by the matrix elements.  These polynomial functions are
nicely behaved with respect to the comultiplication (you get a Hopf
algebra) but we need e.g. continuous functions tending (fast enough)
to zero at infinity in order to use the Haar measure.

\section{Compact Quantum Groups~: Definition and examples}

As we have seen in the previous section, a compact quantum group is a pair
$(A, \comult)$ of a \cst-algebra with an identity and a comultiplication
$\comult$ on $A$ with some extra properties.  We will discuss these extra
conditions on $\comult$ in this section.  We will also show how these
conditions can be verified in examples.

First observe once more that $M(A \ten A)$ now coincides with $A \ten A$
because we have an identity.  So a comultiplication $\comult$ on $A$ is now
a \st-homomorphism from $A$ to $A \ten A$.  Such a \st-homomorphism is
non-degenerate if and only if it is unital.  So in this case, we obtain
the following, easier characterization of a comultiplication.

\blemma
If $A$ is a \cst-algebra with an identity, then a comultiplication on $A$
is a unital \st-homomorphism $\comult$ from $A$ to $A \ten A$ satisfying
coassociativity $(\comult \ten \iota)\comult = (\iota \ten \comult)\comult$.
\elemma

As before, we use the minimal \cst-tensor product.  Again $\iota$ is the
identity map and here $\comult \ten \iota$ and $\iota \ten \comult$ are
the continuous extensions to $A \ten A$ of the obvious maps on the algebraic
tensor product.

Following the ideas of the previous section, a pair $(A, \comult)$ of a
\cst-algebra $A$ with an identity and a comultiplication $\comult$ on $A$
is thought of as a compact quantum semi-group.

Now, what are the extra conditions on $\comult$ in order to have a compact
quantum group~?  To require the existence of a counit and an antipode is not
the best choice as we indicated already in the previous section.  We will
come back to this in section~7.  Instead, the notion of a compact quantum
group, as defined by Woronowicz in~\cite{wor3}, is motivated by the following two
results.

First we have a group-theoretical result.  It can be found in~\cite{hof},
but we
include the simple proof for convenience of the reader.

\bprop \label{prop3.2}
A compact semi-group $G$ with cancellation is a group.
\eprop

\bproof
Take any element $s \in G$. Consider the closed semi-group $H$
generated by $s$.  The intersection of two closed ideals of $H$ is 
again such an ideal, and because $H$ is generated by $s$, 
it is non-empty.  Hence, by
compactness, the intersection $I$ of all closed ideals of $H$
is non-empty.  So it is the smallest closed ideal in $H$.

For any $p \in I$ we have $pI \subseteq I$ and $pI$ is still a closed
ideal in $H$.  Hence $pI=I$ for all $p \in I$.  Given $p \in I$ there is an
element $e \in I$ such that $pe=p$.  Multiply with any $q \in G$ to the right
and cancel $p$ to obtain that $eq=q$.  Then multiply with any element
$p \in G$ to the left and cancel $q$ to get $pe=p$.  So $e$ is the identity
of $G$.

Now, from $se=s$ and the fact that $e \in I$ and that $I$ is an ideal,
we must have $s \in I$.  Because $sI = I$ and $e \in I$, $s$ has an
inverse.  So $G$ is a group.
\eproof

Next, we translate the cancellation law in terms of $\comult$.  We formulate
the result for locally compact semi-groups.

\bprop \label{prop3.3}
Let $G$ be a locally compact semi-group, let $A=C_0(G)$ and $\comult : A \to
M(A \ten A)$ be defined by $(\comult(f))(p,q)=f(pq)$ as before.  Then $G$ has
the cancellation property if the sets $\comult (A)(1 \ten A)$ and $\comult
(A)(A \ten 1)$ are dense subsets of $A \ten A$.
\eprop

\bproof
First let us assume that the set $\comult(A)(1 \ten A)$ is a dense subset of
$A \ten A$.  Take elements $p,q,r$ in $G$ and assume $pr=qr$.  Then, for all
$f,g$ in $A$ we have
$$(\comult(f)(1 \ten g))(p,r) = f(pr)g(r)$$
$$(\comult(f)(1 \ten g))(q,r) = f(qr)g(r)$$
and by assumption, these numbers are equal.  Since this is true for all
$f,g$, again by assumption
$$h(p,r)=h(q,r)$$
for all $h \in A \ten A$.  Therefore $p=q$.  So $G$ has right cancellation.
Similarly, if $\comult(A)(1 \ten A)$ is a dense subset of $A \ten A$, we have
left cancellation.
\eproof

If $G$ is a locally compact group, the maps
$$(p,q) \to (pq,q)$$
$$(p,q) \to (p,pq)$$
are homeomorphisms of $G \times G$.  They give rise to isomorphisms of
$A \ten A$.  These isomorphisms are given by
$$f \ten g \to \comult (f)(1 \ten g)$$
$$f \ten g \to (f \ten 1)\comult(g).$$
So, when $G$ is a locally compact group, the sets $\comult (A)(1 \ten A)$
and $\comult (A)(A \ten 1)$ are dense subspaces of $A \ten A$.

If we combine this result with proposition~\ref{prop3.2}, we find that also the
converse of proposition~\ref{prop3.3} holds for a compact semi-group.

All these observations motivate the following definition of Woronowicz
(given in~\cite{wor3}).

\bdef
A compact quantum group is a pair $(A, \comult)$ of a \cst-algebra $A$ with
identity and a comultiplication $\comult$ on $A$ such that the sets $\comult
(A)(1 \ten A)$ and $\comult (A)(A \ten 1)$ are dense in $A \ten A$.
\end{definition}

From the preceding propositions, we see that any compact group $G$ will give
rise to a compact quantum group and conversely that any compact quantum
group $(A, \comult)$ with $A$ abelian arises in this way.  Recall from the
previous section that this
was one of the properties that we considered desirable for our theory.

Further in the notes, we will show that a nice theory can be developed,
starting from these axioms and that this theory is very similar to the
classical theory.  These are features which we discussed in
the previous section as well,
and which motivate the choice of this definition.

Later in these notes, we will discuss the connection with the definition
of Dijkhuizen and Koornwinder in~\cite{dijk}.  For this we refer to
section~7.

Of course, it is also important that the axioms are relatively easy to
verify in the examples.  This is not so obvious for the density conditions.
Usually such conditions need the existence of some kind of an antipode.
We will show now how this can be done.

We will need the notion of a finite-dimensional representation.

\bdef
Let $(A, \comult)$ be a pair of a \cst-algebra $A$ with identity and a
\st-homo\-morphism $\comult : A \to A \ten A$.  By a finite-dimensional
representation of $(A, \comult)$ we mean a matrix $(v_{pq})$ in $M_n(A)$
for some $n$ such that
$$\comult (v_{pq})= \sum_{k=1}^n v_{pk} \ten v_{kq}$$
for all $p,q=1,\ldots,n$.
If the matrix has an inverse in $M_n(A)$, then the representation is called
non-degenerate.
\end{definition}

We will come back to this definition, and give some motivation in sections
5 and 6 where we study representations more extensively.

Using this notion, we can give the following important result.

\bprop
Let $A$ be a unital \cst- algebra and let $\comult : A \to A \ten A$
be a \st-homo\-morphism.  Assume that $A$ is generated, as a normed algebra,
by the matrix elements of its non-degenerate finite-dimensional
representations. Then $(A, \comult)$ is a compact quantum group.
\eprop

\bproof
Let $(v_{pq})$ be a finite-dimensional representation of $(A, \comult)$.
Then clearly
$$(\iota \ten \comult)\comult(v_{pq}) = (\comult \ten \iota)\comult(v_{pq})
= \sum_{k,l} v_{pk} \ten v_{kl} \ten v_{lq}$$
for all $p,q=1,\ldots,n$.
This implies that $\comult$ is coassociative.

Now assume that $(v_{pq})$ is a non-degenerate finite-dimensional
representation and that $(w_{pq})$ is the inverse of $(v_{pq})$.  Then we
have
$$\sum_k \comult(v_{pk})(1 \ten w_{kq})
= \sum_{k,l}v_{pl} \ten v_{lk}w_{kq}
= \sum_l v_{pl} \ten \delta_{lq}1 = v_{pq} \ten 1.$$
Here $\delta$ denotes the Kronecker delta.
So $v_{pq} \ten 1 \in \comult(A)(1 \ten A)$.
Now, suppose that
\bvgl
a \ten 1 &=& \sum_k \comult(p_k)(1 \ten q_k)\\
b \ten 1 &=& \sum_l \comult(r_l)(1 \ten s_l)
\evgl
with $p_k, q_k, r_l, s_l \in A$.  Then
\bvgl
ab \ten 1 &=& \sum_k \comult (p_k)(1 \ten q_k)(b \ten 1)\\
&=& \sum_k \comult (p_k)(b \ten 1)(1 \ten q_k)\\
&=& \sum_{k,l} \comult (p_k) \comult (r_l) (1 \ten s_l)(1 \ten q_k)\\
&=& \sum_{k,l} \comult (p_k r_l) (1 \ten s_l q_k).
\evgl
So, if $a \ten 1, b \ten 1 \in \comult(A)(1 \ten A)$, then also $ab \ten 1
\in \comult(A)(1 \ten A)$.

Because $v_{pq} \ten 1 \in \comult (A)(1 \ten A)$ when $(v_{pq})$ is a
non-degenerate representation and these matrix elements generate $A$, we
must have that $A \ten 1$ is contained in the closure of $\comult(A)
(1 \ten A)$.  Then of course, $\comult(A)(1 \ten A)$ is dense in $A \ten A$.
Similarly $(A \ten 1)\comult(A)$, and by taking adjoints also $\comult (A)
(A \ten 1)$ is dense in $A \ten A$.
\eproof

Because also the converse of this proposition is true (we will show this in
section~7), one can hope that the result will be useful to construct
examples.

The following proposition is an easy consequence of the previous one.

\bprop \label{prop3.7}
Let $(A, \comult)$ be as in the previous proposition.  Now, assume that $A$
is generated, as a \cst-algebra, by the matrix elements of its
non-degenerate, finite-dimen\-sional representations $(v_{pq})$ so that also
$(v_{pq}^{\ast})$ is an invertible matrix in $M_n(A)$.
Then $(A, \comult)$ is a compact quantum group.
\eprop

\bproof
Remark that
$$\comult(v_{pq}^{\ast}) = \sum_k v_{pk}^{\ast} \ten v_{kq}^{\ast}$$
so that $(v_{pq}^{\ast})$ is also a non-degenerate representation.  Then the
result follows from the previous proposition.
\eproof

We will have a closer look at this representation in section~6.  It will be
called the adjoint representation.  We
will see that it is automatically non-degenerate, if we already have that
$(A, \comult)$ is a compact quantum group.

There is an obvious notion of direct sum of finite-dimensional
representations.  Therefore,
when $A$ is generated by the matrix elements of a finite number of
non-degenerate finite-dimensional representations, it is in fact generated
by the matrix elements of one non-degenerate finite-dimensional
representation.  In this case we have a compact matrix quantum group, as in
\cite{wor1} and~\cite{wor2}. Note that in \cite{wor1} this is called a
compact matrix pseudogroup.

For such a compact matrix quantum group we find the following simple
characterization as a consequence of proposition~\ref{prop3.7}.
It says that a compact matrix quantum group as defined in \cite{wor2}
indeed satisfies the axioms of a compact quantum group given in~\cite{wor3}.

\bprop \label{prop3.8}
Let $A$ be a \cst-algebra with identity and $\comult : A \to A \ten A$ a
\st-homo\-morphism.  Assume that $A$ is generated, as a \cst-algebra, by
elements $(u_{pq})$ such that
$$\comult(u_{pq}) = \sum_k u_{pk} \ten u_{kq}$$
and such that the matrix $(u_{pq})$ and the transpose matrix $(u_{qp})$
are invertible.  Then $(A, \comult)$ is a compact quantum group.
\eprop

\bproof
The transpose matrix $(u_{qp})$ is the conjugate of the matrix
$(u_{pq}^{\ast})$, so the invertibility of the matrix $(u_{qp})$
in $M_n(A)$ is equivalent with the invertibility of $(u_{pq}^{\ast})$.
Now the result follows from proposition~\ref{prop3.7}.
\eproof

Let us now have a look at some examples.
The first example we consider is $SU_{\lambda}(n)$ (see e.g.
\cite{wor.tan} and \cite{vd.dual}).

\bex
Fix $\lambda \in \R$.  Let $A_0$ be the universal unital
\st-algebra generated by
elements $\{u_{pq} \mid p,q=1,...,n \}$ such that
\bvgl
\sum_k u^{\ast}_{kp}u_{kq} &=& \delta_{pq}1\\
\sum_k u_{pk}u^{\ast}_{qk} &=& \delta_{pq}1
\evgl
and such that
$$\sum_{k_i} E(k_1,k_2,...,k_n)u_{l_1k_1}u_{l_2k_2}...u_{l_nk_n}=
E(l_1,l_2,...,l_n)1.$$
Here $E(k_1,k_2,...k_n)=0$ when two indices are equal and otherwise
$$E(k_1,k_2,...k_n) = (-\lambda)^{I(k_1,k_2,...k_n)}$$
where $I(k_1,k_2,...,k_n)$ is the number of inversed pairs in the
permutation $(k_1,k_2,...k_n)$.

We can define a comultiplication $\comult$ on $A_0$ by $\comult(u_{pq})=
\sum u_{pk} \ten u_{kq}$.  It is a straightforward calculation to verify
that the candidates for $\comult(u_{pq})$ satisfy the same relations in
$A_0 \odot A_0$ as the elements $u_{pq}$.  In this way $(A_0,\comult)$
becomes a Hopf \st-algebra.  This was first shown by Woronowicz in
\cite{wor.tan}. However in \cite{vd.dual} a simpler (and more direct)
proof is given.  The counit $\eps$ is given by $\eps(u_{kl})=\delta_{kl}$
and the antipode $\antipode$ by $\antipode(u_{kl})=u_{lk}^{\ast}$.
It is clear that
$\eps$ respects the relations, and that it extends to the counit,
satisfying $(\eps \ten \iota)\comult=(\iota \ten \eps)\comult=\iota$.
The main problem is $\antipode$~: it is not a \st-anti-homomorphism.  But it can
be shown that the adjoints $u_{pq}^{\ast}$ belong to the algebra generated
by the $u_{kl}$, and therefore $\antipode$ is completely determined.  The next
thing then is to show that this candidate for the antipode respects
the defining relations.  The formulas
$$m(\antipode \ten \iota)\comult(u_{pq}) = \eps(u_{pq})1 \mbox{\quad and \quad}
m(\iota \ten \antipode)\comult(u_{pq}) = \eps(u_{pq})1$$
are nothing but the first two defining relations.

In our setting, we do not really need the fact that $(A_0,\comult)$
is a Hopf \st-algebra.  Observe that $(u_{pq})$ is a unitary.
Hence it is invertible in $M_n(A_0)$.  Also (see e.g. lemma 4.7 in
\cite{vd.dual}) the matrix $(u_{pq}^{\ast})$ is invertible.

So we just have to argue that there is a \cst-norm for which
$\comult$ is continuous.
For $a \in A_0$, let
$$\|a\|_0 = \sup \{\|\pi(a)\| \mid
\mbox{$\pi$ is a \st-representation of the algebra $A_0$}\}.$$
Note that the set of \st-representations of $A_0$ is non-empty, take
e.g.\ $\pi(u_{pq})=\delta_{pq}$.
Let $\pi$ be any \st-representation of the
algebra $A_0$ on a Hilbert space $\H$.  Since $(u_{pq})$ is unitary,
we have $\sum_k u_{kq}^{\ast} u_{kq} = 1$.  Hence for all $\xi \in \H$,
$$\|\xi\|^2 = \langle\xi,\xi\rangle = \sum_k
\langle \pi(u_{kq}^{\ast} u_{kq})\xi,\xi\rangle
= \sum_k \langle \pi(u_{kq})\xi,\pi(u_{kq})\xi\rangle
\geq \|\pi(u_{pq})\xi \|^2.$$
So $\|\pi(u_{pq})\| \leq 1$.  This implies that $\|a\|_0$ is finite
for all $a \in A$.  So we have a \cst-seminorm on $A_0$.
Note that $\comult$ is continuous with respect to this seminorm.
Let $I = \{a \in A_0 \mid \|a\|_0 =0\}.$  Then $I$ is a two-sided ideal in
$A_0$.  Define $A$ to be the \cst-completion of $A_0/I$ for the \cst-norm
$\|.\|_0$.  Remark that $\comult$ gives a well-defined comultiplication on
$A_0/I$.  By continuity we can extend this comultiplication to the whole of
$A$.  We will denote this extension again by $\comult$.
Proposition~\ref{prop3.7} now gives that
$(A,\comult)$ is a compact quantum group.
It is a priori not obvious that this quantum group is non-trivial.  However,
it can be shown that there are enough \st-representations of the algebra
$A_0$ (see below).

We will have a closer look at the case $n=2$.  From lemma~4.6 of
\cite{vd.dual} it follows that $u_{12}^{\ast}=- \lambda u_{21}$ and that
$u_{11}^{\ast} = u_{22}$.  If we denote $u_{11}$ by $\alpha$ and $u_{21}$
by $\gamma$, we get
$$u = \left(    \begin{array}{cc}
                        \alpha & -\lambda \gamma^{\ast} \\
                        \gamma & \alpha^{\ast}
                \end{array}   \right)$$
and the defining relations become
$$\alpha^{\ast}\alpha + \gamma^{\ast} \gamma = 1 \hspace{.5in}
\alpha \alpha^{\ast} + \lambda^2 \gamma \gamma^{\ast} =1 $$
$$\alpha \gamma = \lambda \gamma \alpha \hspace{.5in}
\alpha \gamma^{\ast} = \lambda \gamma^{\ast} \alpha \hspace{.5in}
\gamma \gamma^{\ast} = \gamma^{\ast} \gamma. $$
The comultiplication now is given by
\bvgln
\comult(\alpha) &=& \alpha \ten \alpha - \lambda \gamma^{\ast} \ten \gamma\\
\comult(\gamma) &=& \gamma \ten \alpha + \alpha^{\ast} \ten \gamma.
\evgln
Details about the quantum group $SU_{\lambda}(2)$ can be found
in~\cite{wor.SU2}.
It is shown in \cite{wor.SU2} that this \st-algebra can be faithfully
represented by bounded operators on a Hilbert space.  In fact, these
representations of $SU_{\lambda}(2)$ can be used to construct enough
\st-representations of $SU_{\lambda}(n)$ to assure that this quantum
group is non-trivial.
\eex

\bex \label{univ}
This example is the object of study in~\cite{wang.univ}.
Let $Q \in M_n(\C)$ be any invertible matrix.
For a matrix $u=(u_{pq})$ with entries in a $^{\ast}$-algebra, we denote by
$u^t$ the transpose of $u$, by $\overline{u}$ the matrix $(u_{ij}^{\ast})$,
and $u^{\ast}=\overline{u}^t$ the usual conjugation of $u$.
We use the notation $I$ for the $n \times n$ matrix with entries
$\delta_{ij}1$, where $1$ is the identity of the algebra.

Let $A_u(Q)$ be the universal
unital \cst-algebra generated by elements $\{u_{pq} \mid p,q=1,\ldots,n\}$
subject to the following two sets of relations:
$$uu^{\ast} = I = u^{\ast}u$$
and
$$u^t Q \overline{u}Q^{-1} = I = Q \overline{u} Q^{-1} u^t.$$
One can define a comultiplication $\comult$ on $A_u(Q)$ by
$\comult(u_{pq}) = \sum_k u_{pk} \ten u_{kq}$.  It is shown
in~\cite{wang.univ} that the candidates for $\comult(u_{pq})$ satisfy the
same relations in $A_u(Q) \ten A_u(Q)$ as the elements $u_{pq}$.  The second
set of relations for the $u_{pq}$'s implies that the matrix $u^t$ is
invertible.  So by proposition \ref{prop3.8}, $(A_u(Q),\comult)$
is a compact quantum group.

Also here, there is a problem with the \cst-norm.  For some matrices $Q$,
it is in principle possible that the requirement of having a \cst-norm
introduces more relations among the generators so that e.g. in the extreme
case,  the result is a trivial \cst-algebra.  Fortunately, for enough
matrices $Q$, we get non-trivial \cst-algebras (see also \cite{wang.univ}).

These compact matrix quantum groups are universal in the sense that
any compact matrix quantum group is a quantum subgroup of some of these.
Here a compact quantum group $B$ is called a quantum subgroup of another
compact quantum group $A$ if there is a surjective morphism of \cst-algebras
from $A$ to $B$ that preserves the coproducts.
We will give a sketch of the proof of this universality property.
For more details, we refer to~\cite{wang.univ}.

Let $G=(A,\comult)$ be a compact matrix quantum group, generated by the
matrix elements of a non-degenerate finite-dimensional representation $u$.
From theorem~\ref{equnit} it will follow that
we may assume that $u$ is unitary.  This gives
the first set of relations in the definition of the universal compact
matrix quantum groups $A_u(Q)$.

Let $\haar$ be the Haar measure on $(A,\comult)$ (see section 4).
Put $Q=(\iota \ten \haar)
u^t \overline{u}$. Then $Q$ is an invertible matrix in $M_n(A)$
(see e.g. the proof of proposition~\ref{equnit},
in section 6 of this paper, for $v = \overline{u})$.
It is shown in~\cite{wang.univ} that it satisfies
the second set of defining relations in the definition of the
universal compact matrix quantum groups.

\eex

\section{The Haar measure}

A very deep and non-trivial result in the theory of locally compact groups is
the existence and uniqueness (up to a scalar) of the Haar measure.
One of the main difficulties in finding a satisfactory definition for a locally
compact quantum group seems to be exactly this ability to deduce the existence
of (the equivalent of) the Haar measure from the set of axioms.  The fact that
this can be done in the case of compact quantum groups shows again the strength
of Woronowicz' definition.  However, in his proof, Woronowicz assumed the extra
condition that the quantum group admits a faithful state (which is fulfilled in
particular when it is separable).  In~\cite{vd.haar}, it is shown that
 this condition
is not necessary.  We will follow this approach here.  The basic idea of
the proof is the same as in~\cite{wor3}, and is essentially already given
in~\cite{wor1}.  Also the techniques here were already used by Woronowicz.

For a compact group $G$, left invariance of the Haar measure $\haar$ means that for
any $f \in C(G)$
\begin{equation} \label{haar}
\int_G f(ts)\,d \haar(s) = \int_G f(s)\,d \haar(s)
\end{equation}
for all $t \in G$.  This notion can be translated to the corresponding
\cst-algebra $C(G)$.  The Haar measure defines a linear
functional on $C(G)$ which we
will again denote by $\haar$.  Then we have for $F \in C(G \times G)$ that
$$((\iota \ten \haar)F)(t) = \int_G F(t,s)\, d\haar(s)$$
for any $t \in G$.  So equation (\ref{haar})`can be expressed as
$$(\iota \ten \haar) \comult(f) = \haar(f)1$$
for any $f \in C(G)$.  Similarly, right invariance of $h$ will amount
to having
$$(\haar \ten \iota)\comult(f) = \haar(f)1$$
for any $f \in C(G)$.
This justifies the following definition.

\bdef
Let $G=(A, \comult)$ be a compact quantum group.
A linear functional $\haar$ on $A$ is called left invariant if $(\iota \ten
\haar)\comult (a) = \haar(a)1$ for all $a$ in $A$.
A linear functional $\haar'$ on
$A$ is called right invariant if $(\haar' \ten \iota)\comult (a) = \haar'(a)1$
for all $a$ in $A$.
\end{definition}

Let $(A, \comult)$ be a compact quantum group.
We will denote the space of bounded linear functionals on $A$ by $A^{\ast}$.
The comultiplication on $A$ induces a product on $A^{\ast}$ by
$$(\omega_1 \omega_2)(a) = (\omega_1 \ten \omega_2)\comult(a)$$
when $a\in A$ and $\omega_1, \omega_2 \in A^{\ast}$.
We will show that there is a state $\haar$ on $A$ such that $\omega \haar
= h \omega = \omega (1) h$ for any $\omega \in A^{\ast}$.
First, we prove two lemmas.

\blemma \label{haar1}
Let $\omega$ be a state on $A$.  Then there is a state $\haar$ such that
$\haar \omega = \omega \haar = \haar$.
\elemma

\bproof
Define
$$\omega_n = \frac{1}{n}(\omega + \omega^2 + ... + \omega^n).$$
Then $\omega_n$ is a state for all $n$.  Consider a weak \st-limit point
of the sequence $(\omega_n)$.  It is still a state because $A$ has an
identity.  As
$$\|\omega_n \omega - \omega_n\| = \frac{1}{n}\|\omega^{n+1}- \omega\|
\leq \frac{2}{n},$$
we get in the limit $\haar \omega = \haar$.  Similarly $\omega \haar = \haar$.
\eproof

\blemma \label{haar2}
Let $\omega$ and $\haar$ be states on $A$ such that $\haar \omega = \omega
\haar = \haar$.  If $\rho \in A^{\ast}$ and $0 \leq \rho \leq \omega$,
then also $\rho \haar = \rho(1)\haar$.
\elemma

\bproof
Let $a \in A$ and define $b=(\iota \ten \haar)\comult(a)$.  Then
$$(\iota \ten \omega)\comult(b) = (\iota \ten \omega \haar) \comult (a)
= (\iota \ten \haar)\comult (a) =b.$$
So
\bvgl
\lefteqn{(\iota \ten \omega )((\comult (b)-b \ten 1)^{\ast}(\comult (b)-b \ten
1))} \\
\hspace{1cm}&=& (\iota \ten \omega )(\comult(b^{\ast}b))+b^{\ast}b -
(\iota \ten \omega )(\comult(b^{\ast})(b \ten 1)) -
(\iota \ten \omega )((b^{\ast} \ten 1)\comult (b))\\
&=& (\iota \ten \omega )\comult (b^{\ast}b) - b^{\ast}b.
\evgl
If we apply $\haar$ to the above equation, we get
$$ (\haar \ten \omega )((\comult (b)-b \ten 1)^{\ast}(\comult (b)-b \ten 1))=0.$$
Because $0 \leq \rho \leq \omega$, also
$$ (\haar \ten \rho)((\comult (b)-b \ten 1)^{\ast}(\comult (b)-b \ten
1))=0.$$
By the Cauchy-Schwarz inequality, we obtain
$$ (\haar \ten \rho )((c \ten 1)(\comult (b)-b \ten 1))=0$$
for all $c \in A$.  If now we use that $b=(\iota \ten \haar)\comult(a)$,
we get
$$ (\haar \ten \rho \ten \haar)((c \ten 1 \ten 1)(\comult \ten
\iota)\comult(a)) = \rho(1)(\haar \ten \haar)((c \ten 1)\comult(a)). $$
Now
\bvgl
(c \ten 1 \ten 1)(\comult \ten \iota)\comult(a)
&=& (c \ten 1 \ten 1)(\iota \ten \comult)\comult(a)\\
&=& (\iota \ten \comult)((c \ten 1)\comult(a)).
\evgl
By the definition of a compact quantum group, linear combinations of elements
of the form $(c \ten 1)\comult(a)$ are dense in $A \ten A$.  So we can in the
above equality replace $(c \ten 1) \comult(a)$ by $1 \ten q$.  We obtain
$$ (\haar \ten \rho \ten \haar)((\iota \ten
\comult)(1 \ten q)) = \rho(1)(\haar \ten \haar)(1 \ten q), $$
and hence, because $\haar(1)=1$, we get
$$(\rho \ten \haar)(\comult(q)) = \rho(1)\haar(q)$$
for any $q \in A$.
This completes the proof.
\eproof

The existence of an invariant functional now follows easily by using a
compactness argument.

\btheorem
There is a unique left invariant state on $A$.  It is also right invariant.
\etheorem

\bproof
For any $\omega \geq 0$, define
$$K_{\omega} = \{\haar \in A^{\ast} \mid \haar \mbox{ is a state and }
\omega \haar = \omega(1) \haar \}.$$
This is a compact subset of $A^{\ast}$ for the weak \st-topology.  It is
non-empty because of  lemma \ref{haar1}. By lemma \ref{haar2}, we have that
$K_{\omega} \subseteq K_{\rho}$ when $0 \leq \rho \leq \omega$.  It follows
that $K_{\omega_1 + \omega_2} \subseteq K_{\omega_1} \cap K_{\omega_2}$
for any pair of positive linear functionals $\omega_1, \omega_2$.
Therefore the intersection of all $K_{\omega}$ is non-empty.  This proves
the existence of a left invariant state.

Similarly, or by applying this result to the opposite comultiplication
$\comult'$, obtained from $\comult$ by using the flip, there is a state $\haar'$
such that $\haar' \rho = \rho(1)\haar'$ for all $\rho \in A^{\ast}$.  Then
$\haar' \haar = \haar$ by the first formula and $\haar' \haar = \haar'$ by the
second formula.  So $\haar' = \haar$ and the result is proved.
\eproof

So, as in the classical case of a compact group, the left and right
invariant state coincide.  It is called the Haar measure on $(A, \comult)$.

For an example of an explicit calculation of the Haar measure on a
concrete compact quantum group, we refer to~\cite{wor1} and to
\cite{sheu}, where the Haar measures on $SU_{\lambda}(2)$ and on
$SU_{\lambda}(n)$ respectively are calculated.

In general, $\haar$ needs not to be faithful.
The following example gives a compact quantum group with a Haar measure
which is not faithful.

\bex \label{notfaith}
Let $G$ be a discrete group.  Denote by $C^{\ast}(G)$ the full group
\cst-algebra of $G$.  Because $G$ is discrete, this is a unital
\cst-algebra.  For every $s \in G$, we can define $U(s) \in
M(C^{\ast}(G))$ such that
$$(U(s)f)(t) = f(s^{-1}t)$$
for $f \in L^1(G)$.  The mapping $s \mapsto U(s)$ defines a strictly
continuous group representation of $G$ as multipliers on $C^{\ast}(G)$.
Also the mapping $s \mapsto U(s) \ten U(s)$ is a strictly continuous
group representation of $G$ on $M(C^{\ast}(G) \ten C^{\ast}(G))$.
Note that in this context, usually the maximal \cst-norm is used.
In this paper we have chosen to use the minimal \cst-tensor product,
for which the construction of this example remains valid.
It can be shown that one obtains
a unique non-degenerate \cst-homomorphism $\comult : C^{\ast}(G)
\rightarrow M(C^{\ast}(G) \ten C^{\ast}(G))$ such that
$$\comult(U(s)) = U(s) \ten U(s)$$
for all $s \in G$.  Because $(U(s) \ten U(s)) \ten U(s) =
U(s) \ten (U(s) \ten U(s))$, $\comult$ is coassociative.
It is clear that the sets $\comult (C^{\ast}(G))(1 \ten C^{\ast}(G))$ and
$\comult (C^{\ast}(G))(C^{\ast}(G) \ten 1)$ are dense in $C^{\ast}(G)
\ten C^{\ast}(G)$.  Hence $(C^{\ast}(G),\comult)$ is a compact quantum
group.  For the details of this construction, we refer to~\cite{mag}.
In the case that $G$ is a discrete, non-amenable group, the Haar measure
on $(C^{\ast}(G),\comult)$ is not faithful.
\eex

We will come back to the question of faithfulness of the Haar measure in
section~7.

\newpage

\section{The right regular representation}

Let us first define what is meant by a unitary representation of a compact
quantum group $(A, \comult)$.

Let $\H$ be a Hilbert space and denote by $\B_0(\H)$ the \cst-algebra of
all compact operators on $\H$.  Consider the multiplier algebra $M(\B_0(\H)
\ten A)$.  This algebra has two natural embeddings into the multiplier
algebra $M(\B_0(\H) \ten A \ten A)$.  The first one is obtained by extending
the map $x \mapsto x \ten 1$ from $\B_0(\H) \ten A$ to $\B_0(\H) \ten A
\ten A$.  The second one is obtained by composing this map with the
flip on the  last two factors.  It maps $x \ten a$ to $x \ten 1 \ten a$
when $x \in \B_0(\H)$ and $a \in A$.  We will write $v_{(12)}$ and
$v_{(13)}$ for the images of an element $v \in M(\B_0(\H) \ten A)$ by these
two maps.  This is the common {\it leg numbering notation}.
We use brackets to make a distinction with matrix elements.

\bdef
Let $(A, \comult)$ be a compact quantum group.  A representation of $(A,
\comult)$ on a Hilbert space $\H$ is an element $v$ in $M(\B_0(\H) \ten A)$
such that
\begin{equation} \label{rep}
(\iota \ten \comult)(v) = v_{(12)} v_{(13)}.
\end{equation}
If $v$ is unitary, we call it a unitary representation.
\end{definition}

Here $\iota \ten \comult$ is the unique extension to $M(\B_0(\H) \ten A)$
of $\iota \ten \comult$ on the algebraic tensor product.

To see where this definition comes from, consider a strongly continuous
unitary representation $p \mapsto u(p)$ of a compact group on a Hilbert
space $\H$.  Then $u$ is a continuous map from $G$ to $\B(\H)$ where
$\B(\H)$ is considered with the strong operator topology.  Because $u$
is unitary, $u$ is also continuous when we take $\B(\H)$ with the strict
topology, as the multiplier algebra of $\B_0(\H)$.  This means
that $u$ can be considered as an element in $M(\B_0(\H) \ten C(G))$.
We can identify elements in $M(\B_0(\H) \ten C(G) \ten C(G))$ with
strictly continuous $\B(\H)$-valued functions on
$G \times G$.  Then we have
\bvgl
u_{(12)}(p,q) &=& u(p)\\
u_{(13)}(p,q) &=& u(q).
\evgl
Moreover, when $\comult$ is the comultiplication on $C(G)$ given by $(\comult
f)(p,q) = f(pq)$ as before, we also get
$$((\iota \ten \comult)(u))(p,q) = u(pq).$$
So we see that $(\iota \ten \comult)(u) = u_{(12)}u_{(13)}$ means nothing else
but $u(pq) = u(p)u(q)$ for all $p,q \in G$.

A few words about the terminology~:  When the algebra $A$ is non-commutative,
it does not correspond anymore to functions on a compact group.  But we
could think of it as the set of continuous functions on an imaginary
(compact) geometric object.  When we speak of representations of the quantum
group $G=(A, \comult)$, we refer to representations of this underlying
geometric object.  These are in fact the \textit{corepresentations} of the
\textit{coalgebra}.  

Note that, when $v$ is a representation of a compact quantum group
$(A, \comult)$ on a finite-dimensional Hilbert
space $\H$, we can write $v= \sum e_{pq} \ten v_{pq}$ where $(e_{pq})$
are matrix units in $\B(\H)$.  In that case (\ref{rep}) means $\comult(v_{pq})
= \sum_{k=1}^n v_{pk} \ten v_{kq}$.  So we get the notion of
finite-dimensional representation as we defined it in section~3.

Now the Haar measure $\haar$ on the compact quantum group $(A,
\comult)$, obtained in the previous section, will be used to
construct the right regular representation.
This construction and the proofs in this section are essentially copied
from~\cite{wor3}.

Consider the G.N.S. representation of $A$ associated with $\haar$.  We denote the
Hilbert space simply by $\H$, the cyclic vector by $\xi_0$ and we will omit
$\pi_{\haar}$ in the sense that we will write $a \xi$ for $\pi_{\haar}(a) \xi$
when $a \in A$ and $\xi \in \H$.

Let $\K$ be another Hilbert space and assume that $A$ acts faithfully and
non-degenerately on $\K$.  Also here we will write $a \xi$ for the action of an
element $a \in A$ on a vector $\xi \in \K$.

\bprop \label{prop5.1}
There is a unitary operator $u$ on $\H \ten \K$ defined by $u(a \xi_0
\ten \eta) = \comult (a) (\xi_0 \ten \eta)$ when $a \in A$ and $\eta \in \K$.
\eprop

\bproof
If $a_1,a_2,...,a_n \in A$ and $\eta_1,\eta_2,...,\eta_n \in \K$, we have
\bvgl
\| \sum_i \comult(a_i)(\xi_0 \ten \eta_i) \|^2
&=& \sum_{i,j}\langle\comult(a^{\ast}_j a_i)\xi_0 \ten \eta_i, \xi_0 \ten \eta_j\rangle\\
&=& \sum_{i,j}\langle(\haar \ten \iota)\comult(a^{\ast}_j a_i) \eta_i, \eta_j\rangle\\
&=& \sum_{i,j}\haar(a^{\ast}_j a_i) \langle \eta_i, \eta_j\rangle\\
&=& \| \sum_i a_i \xi_0 \ten \eta_i \|^2.
\evgl
It follows that there is an isometric operator $u$ on $\H \ten \K$ given by
$u(a \xi_0 \ten \eta) = \comult (a) (\xi_0 \ten \eta)$ whenever $a \in A$ and
$\eta \in \K$.

Let $a \in A$ and $\eta \in \K$.
From the density of $\comult (A)(1 \ten A)$ in $A \ten A$ it follows
that $a \xi_0 \ten \eta$ can be approximated by linear combinations of
elements $\comult(b)(\xi_0 \ten c\eta)$ with $b,c \in A$.
Now $\comult(b)(\xi_0 \ten c\eta) = u(b \xi_0 \ten c \eta)$.  Since this
holds for any $a \in A$ and $\xi_0$ is cyclic, we see that
the range of $u$ is dense.
Therefore $u$ is unitary.
\eproof

We will now show that $u$ is a multiplier of $\B_0(\H) \ten A$.  Because
$\B_0(\H) \ten A$ acts faithfully and non-degenerately on $\H \ten \K$, and
because finite rank operators are dense in $\B_0(\H)$, it will be sufficient to
show that $u(x \ten 1)$ and $(x \ten 1)u$ are in $\B_0(\H) \ten A$ for all finite
rank operators $x$.

\bprop \label{prop5.2}
The element $u$ is a multiplier of $\B_0(\H) \ten A$.
\eprop

\bproof
Let $x$ be a rank one operator of the form
$x \xi = \langle\xi,\xi_1\rangle a\xi_0$ where
$a \in A$ and $\xi_1,\xi \in \H$.  Then, for $\xi \in \H$
and $\eta \in \K$, we have
\bvgl
u(x \ten 1)(\xi \ten \eta) &=& \langle\xi,\xi_1\rangle u(a\xi_0 \ten \eta)\\
&=&\langle\xi,\xi_1 \rangle \comult (a) \xi_0 \ten \eta.
\evgl
Now, we can approximate $\comult (a)$ by elements of the form $\sum b_i \ten
c_i$.  We have
$$\langle\xi,\xi_1\rangle(\sum b_i \ten c_i)(\xi_0 \ten \eta) =
(\sum y_i \ten c_i)(\xi \ten \eta),$$
where $y_i$ is the rank one operator given by $y_i \xi = \langle\xi,
\xi_1\rangle b_i \xi_0$.
One can check that
$$\left\|u(x \ten 1) - \sum y_i \ten c_i \right\|
\leq \| \xi_1 \| \| \xi_0 \| \left\|\comult(a) - \sum b_i
\ten c_i \right\|.$$
This proves that $u(x \ten 1) \in \B_0(\H) \ten A$.
It follows by the density of $A \xi_0$ in $\H$ that this
will be true for all rank one operators $x$.

Again let $x$ be of the form $x\xi = \langle\xi,  \xi_1\rangle a\xi_0$ and now consider
$u^{\ast}(x \ten 1)$.  For $x \in \H$ and $\eta \in \K$ we get
$$u^{\ast}(x \ten 1)(\xi \ten \eta) = \langle\xi,  \xi_1\rangle u^{\ast} (a \xi_0 \ten \eta).$$
By density of $\comult (A) (1 \ten A)$, we can approximate $a \ten 1$ by
elements of the form $\sum_i\comult(b_i)(1 \ten c_i)$.  Now
\bvgl
\langle\xi,\xi_1 \rangle u^{\ast}  \sum \comult (b_i) (1 \ten c_i) (\xi_0 \ten \eta)
&=& \langle\xi,\xi_1 \rangle \sum b_i \xi_0 \ten c_i \eta\\
&=& \left( \sum y_i \ten c_i \right) (\xi \ten \eta),
\evgl
where $y_i \xi = \langle\xi, \xi_1\rangle b_i \xi_0$.  Again, one can check that
$$\left\|u^{\ast}(x \ten 1) - \sum y_i \ten c_i \right\|
\leq \| \xi_1 \| \| \xi_0 \| \left\|a \ten 1 -
\sum \comult(b_i)  (1 \ten c_i)\right\|.$$
It follows that $u^{\ast}(x \ten 1) \in \B_0(\H) \ten A$ and again
by density, that this holds for all rank one operators $x$.

By the remark preceding the proposition, this completes the proof.
\eproof

Now, we can consider $(\iota \ten \comult)(u)$ and prove that $u$ is a
representation.

\btheorem
We have $(\iota \ten \comult)(u) = u_{(12)} u_{(13)}$ and therefore, $u$
is a unitary representation of the compact quantum group $(A, \comult)$.
\etheorem

\bproof
Let $b,c \in A$ and $\eta_1, \eta_2 \in \K$.  Then
\bvgl
u_{(12)} (b \ten 1 \ten c)(\xi_0 \ten \eta_1 \ten \eta_2) &=&
\comult(b) (\xi_0 \ten \eta_1) \ten c \eta_2\\
&=& (\comult \ten \iota)(b \ten c)(\xi_0 \ten \eta_1 \ten \eta_2).
\evgl
Let $a \in A$.
Again because we can approximate $\comult(a)$ by elements of the form
$\sum b_i \ten c_i$, we can replace $b \ten c$ in this equation by
$\comult(a)$.  We get
$$u_{(12)}u_{(13)}(a \xi_0 \ten \eta_1 \ten \eta_2)
= ((\comult \ten \iota)\comult(a))(\xi_0 \ten \eta_1 \ten \eta_2).$$
On the other hand, we have for any $y \in \B(\H)$ and $b \in A$ that
$$(\iota \ten \comult)(y \ten b)(a \xi_0 \ten \eta_1 \ten \eta_2)
= ya \xi_0 \ten \comult (b) (\eta_1 \ten \eta_2).$$
Let $x \in \B_0(\H)$ be any compact operator on $\H$.
Then $(x \ten 1)u \in \B_0(\H)$ can be approximated by linear combinations
of elements $y \ten b$.  Hence we can replace in the previous
formula $y \ten b$ by $(x \ten 1)u$.  We get
$$(x \ten 1 \ten 1)((\iota \ten \comult)(u))(a \xi_0 \ten \eta_1 \ten \eta_2)
= (x \ten 1 \ten 1)((\iota \ten \comult) \comult(a))
(\xi_0 \ten \eta_1 \ten \eta_2).$$
Because $x$ was arbitrary, we can cancel it.
By the coassociativity of $\comult$, we get $(\iota \ten \comult)(u)
=u_{(12)} u_{(13)}$ and this proves the result.
\eproof

It is not hard to see that in the group case, we indeed get the right
regular representation here.  If $G$ is a compact group, let $A=C(G)$
 and $\comult$ be as before.  In this case, $\H = L^2(G)$, the Hilbert
space of (classes of) square integrable functions on $G$ with
respect to the Haar measure, and $C(G)$ acts on $\H$ by multiplication.
Classically, the right regular representation $u: G \rightarrow \B(\H)$ of
$G$ is given by
$$ (u(q)\xi)(p) = \xi (pq)$$
when $p,q \in G$ and $\xi \in L^2(G)$.  We will consider $u$ as an
element in $M(\B_0(\H) \ten C(G))$, and let this algebra act on $\H
 \ten \H$.  Then we have for $\xi, \eta \in \H$ and $p,q \in G$
that
\bvgl
u(\xi \ten \eta)(p,q) &=& \xi(pq) \eta(q)\\
&=& (\comult (\xi)(1 \ten \eta))(p,q).
\evgl
Remembering that here the cyclic vector $\xi_0 =1$, we see that the
notion of right regular representation for a compact quantum group indeed
corresponds to the classical one in the case of a compact group.

For the explicit calculation of the regular representation in an example,
we refer to \cite{lan}, where Lance constructed the regular
representation of $SU_{\lambda}(2)$.

Next, we prove a technical result which will be very useful in the
following sections.

\bprop  \label{prop.dense}
The set $\{(\omega \ten \iota)(u) \mid \omega \in \B_0(\H)^{\ast} \}$ is
dense in $A$.
\eprop

\bproof
Let $a,b \in A$.  Recall that for given vectors $\xi_1,\xi_2 \in \H$, we denote
by $\omega_{\xi_1,\xi_2}$ the linear functional on $\B(\H)$ such that
$\omega_{\xi_1,\xi_2}(x)=\langle x\xi_1,\xi_2 \rangle$ for all $x \in \B(\H)$.
We have
\bvgl
\langle((\omega_{a \xi_0, b \xi_0} \ten \iota)u) \eta_1, \eta_2\rangle
&=& \langle u (a \xi_0 \ten \eta_1), b \xi_0 \ten \eta_2\rangle\\
&=& \langle\comult(a) (\xi_0 \ten \eta_1), b \xi_0 \ten \eta_2\rangle\\
&=& \langle(\haar \ten \iota)((b^{\ast} \ten 1)\comult (a))\eta_1, \eta_2\rangle
\evgl
for any $\eta_1, \eta_2 \in \K$.  So
$$ (\omega_{a \xi_0, b \xi_0} \ten \iota)u =
(\haar \ten \iota)((b^{\ast} \ten 1)\comult (a)).$$
By the axioms of a compact quantum group, linear combinations
of elements of the form $(b^{\ast} \ten 1) \comult(a)$ with $a,b \in A$
are dense in $A \ten A$.  It follows that $\{(\omega_{a\xi_0, b\xi_0}
\ten \iota)(u) \mid a,b \in A \}$ is dense in $A$.
Remark also that by density of $A\xi_0$ in $\H$, it follows that
$(\omega \ten \iota)(u) \in A$ for any $\omega$ of the form
$\omega_{\xi_1,\xi_2}$ with $\xi_1,\xi_2 \in \H$.  And because linear
combinations of these functionals are dense in $\B_0(\H)^{\ast}$, we get that
$(\omega \ten \iota)(u) \in A$ for any $\omega \in \B_0(\H)^{\ast}$.
\eproof

Now we show that $u$ implements $\comult$.

\bprop
For all $a \in A$ we have that $\comult (a) = u(a \ten 1)u^{\ast}$.
\eprop

\bproof
When $a,b \in A$ and $\eta \in \K$ we have
\bvgl
u(a \ten 1)(b \xi_0 \ten \eta) &=& \comult(ab) (\xi_0 \ten \eta)\\
&=& \comult(a) \comult(b) (\xi_0 \ten \eta)\\
&=& \comult(a) u(b \xi_0 \ten \eta).
\evgl
Hence $u(a \ten 1) = \comult(a) u$ and because $u$ is unitary,
$\comult (a) = u (a \ten 1)u^{\ast}$.
\eproof

This result remains true if we take $\K = \H$.  And this will imply that
the element $u$ is a multiplicative unitary in the sense of~\cite{baaj}.
We can deduce for
$u$, which is in $M(\B_0(\H) \ten A)$, that
$$(\iota \ten \comult)(u) = (1 \ten u)(u \ten 1)(1 \ten u^{\ast}).$$
On the other hand we know that $(\iota \ten \comult)(u) = u_{(12)}u_{(13)}$.
So we have
$$(1 \ten u)(u \ten 1)(1 \ten u^{\ast}) = u_{(12)}u_{(13)},$$
or equivalently,
$$u_{(23)}u_{(12)} = u_{(12)}u_{(13)}u_{(23)}.$$
So $u$ is a multiplicative unitary.

Using the arguments of the proofs of propositions \ref{prop5.1} and
\ref{prop5.2} one can
show that there exists a unitary $v \in M(\B_0(\H) \ten A)$ such that
$$v^{\ast}(a \xi_0 \ten \eta) = \comult' (a) (\xi_0 \ten \eta)$$
when $a \in A$ and $\eta \in \K$.  Here $\comult'$
denotes the opposite comultiplication $\sigma \circ \comult$, where
$\sigma$ is the flip map.
Again, $v$ is a representation
and will be called the left regular representation of $(A, \comult)$.
Also here, it is not hard to see that in the group case, we indeed
get the left regular representation.

\section{Representation theory}

\setcounter{equation}{0}

In the previous section we have constructed the right regular representation
by using the Haar measure.  The main result in this section is that any
unitary representation can be decomposed into irreducible unitary
representations and that all irreducible unitary representations
are contained in the regular representation.

Let $(A,\Phi)$ be a compact quantum group.  Let $v$ be a representation of
$(A,\Phi)$ on a Hilbert space $\H$.  A closed subspace $\H_1$ of $\H$ is said
to be invariant if $(e \ten 1)u(e \ten 1) = u(e \ten 1)$, where $e$ is the
orthogonal projection onto this subspace.  The representation $v$ is
called irreducible if the only invariant subspaces are $\{0\}$ and $\H$.
Like in the classical case, also here, when $u$ is a unitary representation
and $\H_1$ is an invariant subspace, also the orthogonal complement of $\H_1$
will be invariant.  We will first obtain this result.  Throughout this
section, we will denote the Haar measure on $(A, \comult)$ by $\haar$.

\bprop \label{propos}
Let $v$ be a unitary representation of $(A,\Phi)$ on a Hilbert space $\H$.
Denote by $\B$ the normclosure in $\B(\H)$ of the set of elements $(\iota
\ten h)(v(1 \ten a))$ where $a \in A$.  Then $\B$ is a non-degenerate
\cst-subalgebra of $\B(\H)$.  Moreover $v \in M(\B \ten A)$.
\eprop

\bproof
Take $a \in A$ and define an element $x \in \B(\H)$ by $x=(\iota\ten h)(v(1
\ten a))$.  Using the left invariance of $h$ and the fact that $v$ is a
representation, one can easily show that
\begin{equation} \label{hulp}
v^\ast (x \ten 1) = (\iota \ten h)(v_{(13)}(1 \ten \Phi(a)))
\end{equation}
in $M(\B_0(\H) \ten A)$.

When $a,b \in A$ and $x,y$ are defined by
\bvgl
x &=& (\iota \ten h)(v (1 \ten a))\\
y &=& (\iota \ten h)(v (1 \ten b)),
\evgl
using result (\ref{hulp}), we get
$$ y^\ast x = (\iota \ten h \ten h)(v_{(13)} (1 \ten (b \ten 1)\Phi(a))).$$
It follows that $y^\ast x = (\iota \ten h)(v(1 \ten c))$ where
$c = (h \ten \iota)((b \ten 1)\Phi(a))$.
This implies that $\B^\ast \B \subseteq \B$.  But as $(A \ten 1)\Phi(A)$ is
dense in $A \ten A$, it also implies that $\B^\ast \B$ is dense in $\B$.  As a
consequence, $\B$ is self-adjoint and $\B \B \subseteq \B$.  Hence $\B$ is a
\cst-subalgebra of $\B(\H)$.

To prove that $\B$ acts non-degenerately on $\H$, take any $x \in \B_0(\H)$
and $a \in A$.  Then $(\iota \ten h)(v(x \ten a)) \subseteq \B \B_0(\H)$.
Because $v$ is a unitary multiplier of $M(\B_0(\H) \ten A)$, we have that
$v(\B_0(\H) \ten A)$ is dense in $\B_0(\H) \ten A$.  It follows that $\B
\B_0(\H)$ is dense in $\B_0(\H)$.  This proves that $\B$ acts
non-degenerately.

Finally, we show that $v \in M(\B \ten A)$.  Again using the
formula~(\ref{hulp}), we find for every $a,b \in A$ that
$$v^\ast(x \ten b) = (\iota \ten \iota \ten h)(v_{(13)}(1 \ten (\Phi(a)(b
\ten 1))))$$
when $x = (\iota \ten h)(v(1 \ten a))$.  The right hand side belongs to $\B
\ten A$ as $\Phi(a)(b \ten 1) \in A \ten A$.  Moreover, as $\Phi(A)(A \ten
1)$ is dense in $A \ten A$, we have that $v^\ast(\B \ten A)$ is dense in
$\B \ten A$.  Then also $v (\B \ten A) \subseteq \B \ten A$ and we obtain
that $v \in M(\B \ten A)$.
\eproof

Now we can prove that also orthogonal complements of invariant subspaces
are invariant.

\bprop
Let $v$ be a unitary representation of $(A,\Phi)$ on $\H$.  If $\H_1$ is an
invariant subspace, then also the orthogonal complement of $\H_1$ is
invariant.
\eprop

\bproof
Denote the orthogonal projecton of $\H$ onto $\H_1$ by $e$.  Since $v(e
\ten 1) = (e \ten 1)v(e \ten 1)$, we will have $xe = exe$ for every $x \in
\B$, where $\B$ is the algebra defined in the previous proposition.
As $\B$ is self-adjoint, this implies $xe = ex$ for all $x \in \B$.
As $\B$ acts non-degenerately on $\H$, it follows that $xe = ex$ for all $x
\in M(\B)$.  Then also $v (e \ten 1) = (e \ten 1)v$.
\eproof

The theory of representations of a compact quantum group is very similar to
the one of compact groups.  The result above illustrates this.  In the
classical situation, the proofs are of course easier.  But except for the
result above, in most cases, the general theory is not much more
complicated than the classical one.

Just as for compact groups, also for compact quantum groups, it is possible
to define the direct sum and the tensor product of two representations. The
definition for the direct sum is more or less obvious.
The tensor product of two representations $v$ and $w$
of $(A, \comult)$ on Hilbert spaces $\H_1$ and $\H_2$ respectively, is
defined as
$v \tenrep w = v_{(13)}w_{(23)}$.
Then $v \tenrep w \in M((\B_0(\H_1) \ten \B_0(\H_2)) \ten A)
\cong M(\B_0(\H_1 \ten \H_2) \ten A)$, and one can check that it is again
a representation.

Let $v,w$ be representations
of $(A, \comult)$ on Hilbert spaces $\H_1$ and $\H_2$ respectively.
An intertwiner between $v$ and $w$ is an element $x \in \B(\H_1,\H_2)$
such that $(x \ten 1)v=w(x \ten 1)$.
We will denote the set of intertwiners between $v$ and $w$ by
$\mbox{Mor}(v,w)$.
Two representations are said to be equivalent
if there is an invertible intertwiner.  They are unitarily equivalent if
the intertwiner can be chosen to be unitary.

We will use the Haar measure, as in the classical  theory, to
construct intertwiners.

\blemma \label{lemma6.1}
Let $v$ and $w$ be representations of $(A, \comult)$ on Hilbert spaces
$\H_1$ and $\H_2$ respectively.  Let $x$ be any compact operator from
$\H_1$ to $\H_2$.  Define
$$y = (\iota \ten \haar)(w^{\ast}(x \ten 1)v).$$
Then $y$ is a compact operator from $\H_1$ to $\H_2$ such that
$$w^{\ast}(y \ten 1)v = y \ten 1.$$
\elemma

\bproof
Denote by $\B_0(\H_1, \H_2)$ the compact operators from $\H_1$ to $\H_2$.
Because $x \in \B_0(\H_1,\H_2)$ and $v \in M(\B_0(\H_1) \ten A)$ we have
$(x \ten 1)v \in \B_0(\H_1,\H_2) \ten A$ and also \newline
$w^{\ast}(x \ten 1)
v \in \B_0(\H_1,\H_2) \ten A$.  Therefore $y$ is well defined and $y \in
\B_0(\H_1, \H_2)$.

We also have
$$(\iota \ten \comult)(w^{\ast}(x \ten 1)v) = w^{\ast}_{(13)}
w^{\ast}_{(12)}(x \ten 1 \ten 1)v_{(12)}v_{(13)}.$$
If we apply $\iota \ten \haar \ten \iota$ to the left hand side of this
equation we get $y \ten 1$ by the invariance of $\haar$.  If we apply
the same map to the right hand side, we obtain $w^{\ast}(y \ten 1)v$.
This proves the lemma.
\eproof

If $w$ is a unitary representation,
then $(y \ten 1)v = w (y \ten 1)$ and $y$ will
be an intertwiner between $v$ and $w$.

\bprop \label{equnit}
Any non-degenerate finite-dimensional representation is equivalent with
a unitary representation.
\eprop

\bproof
So, let $\H$ be a finite-dimensional Hilbert space and assume that $v$ is
a non-degenerate representation on $\H$. Let
$$y = (\iota \ten \haar)(v^{\ast}v).$$
By the previous lemma, $v^{\ast}(y \ten 1)v = y \ten 1$.
Now $y \geq 0$ and because $v$ is
invertible, $v^{\ast}v \geq \eps 1$ for some $\eps > 0$ and so
$y \geq \eps 1$.  This implies that $y$ is invertible.
If we let
$$w = (y^{1/2} \ten 1) v (y^{-1/2} \ten 1)$$
we get a representation, equivalent with $v$.
Because $w^{\ast}w = 1$ and $w$ is invertible, we have that $w$ is unitary.
This proves the proposition.
\eproof

The result is still true for infinite-dimensional representations.  Then
we must use an extended version of lemma~\ref{lemma6.1}
with $x$ a bounded operator.

Next we show that any unitary representation decomposes into irreducible
finite-dimensional representations.

\btheorem \label{th6.3}
Let $v$ be a unitary representation of $(A, \comult)$ on a Hilbert space $\H$.
Then there is a set $\{e_{\alpha} \mid \alpha \in I\}$ of mutually orthogonal
finite-dimensional projections with sum $1$ and satisfying
$$v(e_{\alpha} \ten 1) = (e_{\alpha} \ten 1)v$$
and $v(e_{\alpha} \ten 1)$, considered as an element in $\B(e_{\alpha} \H)
\ten A$ is a finite-dimensional irreducible unitary representation of
$(A, \comult)$.
\etheorem

\bproof
Consider the set $\B$ of operators $y$ in $\B_0(\H)$ satisfying $v(y \ten 1)
= (y \ten 1)v$.  This is clearly a norm closed subalgebra of $\B_0(\H)$.
Because $v$ is unitary, it is also a self-adjoint subset of $\B_0(\H)$.
Therefore $\B$ is a \cst-subalgebra of $\B_0(\H)$.

By lemma~\ref{lemma6.1} we see that
$$(\iota \ten \haar)(v^{\ast}(x \ten 1)v) \in \B$$
for all $x \in \B_0(\H)$.  If we consider an increasing net $(x_{\lambda})$
of positive elements in $\B_0(\H)$ with $\sup x_{\lambda} = 1$ then we get
an increasing net $(y_{\lambda})$ of positive elements in $\B$ such that
$$\sup y_{\lambda} = (\iota \ten \haar)(v^{\ast}v) = \haar(1)1 =1.$$
It follows that $\B$ acts non-degenerately on $\H$.  Then, taking a maximal
family of mutually orthogonal, minimal projections in $\B$, we get the
desired result.
\eproof

So any irreducible unitary representation of a compact quantum group is
finite-dimen\-sional.  Now we can easily prove Schur's lemma:

\blemma
Let $v,w$ be irreducible unitary representations of $(A, \comult)$ on Hilbert
spaces $\H_1$ and $\H_2$ respectively.  Then either $v$ is not
equivalent to $w$ and $\mbox{\rm Mor}(v,w)=0$, or $v$ is equivalent to $w$
and there is
an invertible $x \in \B(\H_1,\H_2)$ such that $\mbox{\rm Mor}(v,w)
= \{\lambda x \mid \lambda
\in \C\}$.  In particular, if $u$ is an irreducible unitary representation of
$G$, then $\mbox{\rm Mor}(u,u) = \{\lambda I \mid \lambda \in \C\}$.
\elemma

\bproof
Notice that $\H_1$ and $\H_2$ are finite-dimensional.
It is straightforward to check that the kernel of an intertwiner between
$v$ and $w$ is a $v$-invariant subspace, and that the image of an intertwiner
is a $w$-invariant subspace.

Suppose there is an $x \in \mbox{Mor}(v,w)$ with $x \not= 0$. So the kernel of $x$ is
a $v$-invariant subspace, which is not the whole of $\H$.  Because $v$ is
irreducible, it follows that $\mbox{ker}(x)=\{0\}$.  Similarly, the image
of $x$ is a non-zero $w$-invariant subspace and hence must be the whole of
$\K$.  So $x$ is invertible.

Now let $y \in \mbox{Mor}(v,w)$.  For any $\lambda \in \C$, it is clear that
$\lambda x -y$
is again an intertwiner between $v$ and $w$.  Take $\lambda$ such that
$\mbox{det}(\lambda x -y) = 0$.  Then the kernel of $\lambda x - y$ is
non-zero, and hence equals $\H$.  So $y= \lambda x$.
\eproof

We can apply the result of theorem~\ref{th6.3} to the right regular
representation,
constructed in the previous section.  The next theorem says that all
irreducible unitary representations will be obtained.
Until now the proofs of this section are the ones given in~\cite{wor3}.
The next results were also present in~\cite{wor3}, but the proofs given here
are different (and simpler), using the technique of lemma~\ref{lemma6.1}.

\btheorem \label{th6.4}
Every irreducible unitary representation is contained in the regular
representation. \etheorem

\bproof
Let $u$ be the right regular representation, acting on $\H$ and let $v$ be
any irreducible unitary representation, acting on $\K$.  Take any compact
operator $x$ from $\H$ to $\K$ and consider, as before
$$y = (\iota \ten \haar)(v^{\ast}(x \ten 1)u).$$
By lemma~\ref{lemma6.1} this is an operator from $\H$ to $\K$ such that
$$(y \ten 1)u = v(y \ten 1).$$
Let $e$ be the projection on the range of $y$.  Then
$$(e \ten 1)v(y \ten 1) = (ey \ten 1)u = (y \ten 1)u = v(y \ten 1)$$
so that
$$(e \ten 1)v(e \ten 1) = v(e \ten 1).$$
Because $v$ is assumed to be irreducible, we have $e=0$ or $e=1$.  This
means that $y=0$ or $y$ is surjective.

If $y$ is surjective, it follows that $v$ is equivalent with a
subrepresentation of $u$.
Now, suppose that $y=0$ for all $x$.
Consider any $x$ of the form $x \xi = \langle \xi, \xi_1 \rangle \eta_1$
where $\xi_1 \in \H$ and $\eta_1 \in \K$.
Note that for any $\xi \in \H$ and $\eta \in \K$
$$(v^{\ast}(x \ten 1)u)(\xi \ten \eta) = (v^{\ast}(1 \ten a))
(\eta_1 \ten \eta)$$
with $a=(\omega_{\xi,\xi_1} \ten \iota)u$.  So $y=0$ implies
$(\iota \ten \haar)(v^{\ast}(1 \ten a))(\eta_1)=0$.  Because this holds
for any $\eta_1$ we have
$$(\iota \ten \haar)(v^{\ast}(1 \ten a))=0$$
for any $a$ of the form $(\omega_{\xi, \xi_1} \ten \iota)u$.
Now we have seen in proposition~\ref{prop.dense}
that such elements are dense in $A$.
Therefore
$$(\iota \ten \haar)(v^{\ast}(1 \ten a)) = 0$$
for all $a \in A$.  If we multiply to the right with any $x \in \B(\K)$ and
use linearity we find that also $(\iota \ten \haar)(v^{\ast}v) = 0$.
But this is a contradiction as $v^{\ast}v = 1$.
\eproof

Theorem \ref{th6.3} and \ref{th6.4} give the main results in this section.
However, in the
next section, we will also need the adjoint of a unitary representation.

\bdef
Let $v$ be a representation on a finite-dimensional Hilbert space $\H$.
If we denote the matrix units in $\B(\H)$ by $(e_{pq})$,
we can write $v = \sum e_{pq}
\ten v_{pq}$.  Define $\overline{v} = \sum e_{pq} \ten v_{pq}^{\ast}$.
Then $\overline{v}$ is still a representation. It is called the adjoint
of $v$.
\end{definition}

The definition of $\overline{v}$ depends on the choice of the matrix units.
But a different choice will give an equivalent representation.  So, the
adjoint is only defined up to equivalence.

In the classical theory, the adjoint of a unitary representation is again
a unitary representation.  In the quantum case, this needs no longer be
true.  We will have however that it is equivalent to a unitary and that
it is still irreducible.  To obtain this result, we would need e.g.\ that
$\overline{v}$ is invertible when $v$ is unitary.  This seems to be not
so obvious.  We have to give a proof using the basic technique of lemma
~\ref{lemma6.1} again.

\blemma
If $v$ is a finite-dimensional irreducible representation, then
$\overline{v}$ is also irreducible.
\elemma

\bproof
Let $v$ act on $\H$.  Choose matrix units $(e_{pq})$ in $\B(\H)$ and define
the anti-homomorphism $\gamma$ of $B(\H)$ by $\gamma(e_{pq}) = e_{qp}$.
Remark that $\overline{v} = ( \gamma \ten \iota) v^{\ast}$.  Now let $e$
be a projection in $B(\H)$ such that
$$\overline{v}(e \ten 1) = (e \ten 1) \overline{v} (e \ten 1).$$
If we apply $\gamma \ten \iota$, we get
$$(\gamma (e) \ten 1) v^{\ast} = (\gamma(e) \ten 1) v^{\ast}
(\gamma(e) \ten 1)$$
and if we take adjoints, we obtain
$$v(f \ten 1) = (f \ten 1)v(f \ten 1)$$
where $f = \gamma(e)^{\ast}.$  Then $f$ is also a projection and because
$v$ is irreducible, it is either $0$ or $1$.  This implies that also $e=0$
or $e=1$.
\eproof

\bprop
If $v$ is an irreducible unitary representation, then $\overline{v}$ is
equivalent with a unitary representation.
\eprop

\bproof
Let $u$ be the right regular representation, acting on $\H$.  Let $v$
act on the finite-dimensional space $\K$.  Take any compact operator
$x$ from $\K$ to $\H$ and define as before
$$y = (\iota \ten \haar)(u^{\ast}(x \ten 1)\overline{v}),$$
where $\haar$ denotes again the Haar measure.
This is an operator from $\K$ to $\H$ such that $(y \ten 1)\overline{v}
= u(y \ten 1)$.  We will show that $y \not= 0$ for some $x$.  Because
$\overline{v}$ is irreducible, this will imply that $\overline{v}$ is
equivalent with a subrepresentation of $u$.

So assume that $y=0$ for all $x$.  This means
$$(\iota \ten \haar)(u^{\ast}(x \ten 1) \overline{v}) = 0$$
for all operators $x$ from $\K$ to $\H$.
As in the proof of theorem~\ref{th6.4}, this implies
$$(\iota \ten \haar)(u^{\ast}(1 \ten b)) = 0$$
where $b=(\psi \ten \iota)\overline{v}$ for any $\psi \in
\B(\K)^{\ast}$ (recall that $\K$ is finite-dimensional, so any operator on
$\K$ has finite rank).
By applying any $\omega \in \B_0(\H)^{\ast}$, we get
$$\haar(ab) = 0$$
where $a=(\omega \ten \iota)u^{\ast}$.  By lemma~\ref{prop.dense},
elements $a$ of this form are dense in
$A$.  So we obtain
$$(\psi \ten \iota)(\iota \ten \haar)((1 \ten a) \overline{v}) = 0$$
for any $\psi \in \B(\K)^{\ast}$.  Hence
$$(\iota \ten \haar)((1 \ten a) \overline{v}) = 0$$
for all $a \in A$.  By applying transposition $\gamma$ we get
$$(\iota \ten \haar)((1 \ten a) v^{\ast}) = 0$$
for all $a \in A$.  Then also
$$(\iota \ten \haar)((x \ten a)v^{\ast}) = 0$$
for all $x \in \B(\K)$ and $a \in A$, and by linearity
$$(\iota \ten \haar)(vv^{\ast}) = 0.$$
This is a contradiction as $vv^{\ast} =1$.
\eproof

Remark that it follows that $\overline{v}$ is non-degenerate whenever $v$
is non-degenerate.

\section{The Hopf $^{\ast}$-algebra of matrix elements.}

Let $A_0$ be the subspace of $A$ spanned by the matrix elements of
finite-dimensional unitary representations.  We will show in this section
that $A_0$ is a dense $^{\ast}$-subalgebra of $A$, that $\comult(A_0)$ is
contained in the algebraic tensor product
$A_0 \odot A_0$ and that $(A_0, \comult)$ is a Hopf $^{\ast}$-algebra.
This result comes from~\cite{wor3}, but the proofs given here are
somewhat different.

In~\cite{dijk},
Dijkhuizen and Koornwinder start from this Hopf $^{\ast}$-algebra $A_0$
to define a compact quantum group. We will have a look at their
definition in remark~\ref{dijkhuizen}.
In this section we will also show the
equivalence between Woronowicz' two definitions of a compact matrix
quantum group.  We end the section with a discussion about the
faithfulness of the Haar measure.

\bprop
$A_0$ is a dense $^{\ast}$-subalgebra of $A$.
\eprop

\bproof
The product of two matrix elements of finite-dimensional unitary representations
is a matrix element of the tensor product of these two representations.
Therefore $A_0$ is a subalgebra.

Because the adjoint of a finite-dimensional unitary representation is equivalent
with a unitary representation, $A_0$ is $^{\ast}$-invariant.

Let $u$ be the right regular representation, acting on the Hilbert space $\H$.
From theorem \ref{th6.3} it follows that $u$ decomposes into irreducible
finite-dimensional representations $\{u_{\alpha} \mid \alpha \in I\}$ acting on
subspaces $\{\H_{\alpha} \mid \alpha \in I \}$.  For $\alpha \in I$, let
$n(\alpha)$ be the dimension of $\H_{\alpha}$ and let
$\{ \xi_1^{\alpha},\, \xi_2^{\alpha},\,...\, ,\, \xi_{n(\alpha)}^{\alpha} \}$
be an
orthonormal basis  for $\H_{\alpha}$.  Define $\omega^{\alpha \beta}_{pq}
\in \B_0(\H)^{\ast}$ by $\omega^{\alpha \beta}_{pq}(x) = \langle
x \xi^{\alpha}_p, \xi^{\beta}_q \rangle$ for $x \in \B_0(\H)$.
Then the linear span of
$$\{ \omega^{\alpha \beta}_{pq} \mid
\alpha, \beta \in I,\ 1 \leq p \leq n(\alpha), \ 1 \leq q \leq n(\beta) \}$$
is dense in $\B_0(\H)^{\ast}$.  Hence it follows from proposition
\ref{prop.dense} that the linear span of the set
$$\{ (\omega^{\alpha \beta}_{pq} \ten \iota) u \mid
\alpha, \beta \in I,\ 1 \leq p \leq n(\alpha),\ 1 \leq q \leq n(\beta) \}$$
is dense in $A$.  Now $(\omega^{\alpha \beta}_{pq} \ten \iota) u = 0$
unless $\alpha = \beta$ and the elements
$(\omega^{\alpha \alpha}_{pq} \ten \iota) u$
are exactly the matrix elements of the representation $u_{\alpha}$.  This
proves the proposition.
\eproof

Remark that $1 \in A_0$ because $1$ is a representation.

Take a complete set $\{u^{\alpha} \mid \alpha \in I \}$ of mutually inequivalent,
irreducible unitary representations.  We want to show that the
elements
$\{ u^{\alpha}_{pq} \mid \alpha \in I,\ 1 \leq p,q \leq n(\alpha) \}$
form a basis for $A_0$.  In order to prove the linear independence, we first
deduce the following orthogonality relations.

\bprop \label{orth}
Let $\haar$ be the Haar
measure on $(A, \comult)$.  For any $\alpha \in I$, there is a positive
invertible operator $F^{\alpha}$ acting on $\H_{\alpha}$ such that for any
$\alpha, \beta \in I$ and $1 \leq j,q \leq n(\alpha)$,
$1 \leq i,p \leq n(\beta)$
$$\haar((u^{\beta}_{ip})^{\ast}u^{\alpha}_{jq}) = \delta_{\alpha \beta}
\delta_{pq} F^{\alpha}_{ij}.$$
\eprop

\bproof
We apply lemma~\ref{lemma6.1} for $u^{\alpha}$ and $u^{\beta}$, with
$x$ the operator mapping the $j$th basis element of $\H_{\alpha}$ to the
$i$th basis element of $\H_{\beta}$ and the other basis elements to zero.
We obtain that the matrix with elements
$$y_{pq} = \haar ((u^{\beta}_{ip})^{\ast} u^{\alpha}_{jq})$$
intertwines $u^{\alpha}$ with $u^{\beta}$.
Because $u^{\alpha}$ and $u^{\beta}$ are irreducible we have
$y_{pq} = 0$ if $\alpha \not= \beta$.  In the case $\alpha = \beta$ there is
$\lambda \in \C$ such that $y_{pq} = \lambda \delta_{pq}$
for all $p,q$.  This $\lambda$ depends on $\alpha$ and $i,j$.  Let us
write
$$F^{\alpha}_{ij} \delta_{pq} = \haar ((u^{\alpha}_{ip})^{\ast}
u^{\alpha}_{jq}).$$
Now we apply lemma~\ref{lemma6.1} for $\overline{u^{\alpha}}^{\ast}$ and
$(u^{\alpha})^t$.
We obtain that $(\iota \ten \haar)(\overline{u^{\alpha}}(u^{\alpha})^t)$
is an intertwiner
between $\overline{u^{\alpha}}$ and $(u^{\alpha})^t$.
Because $(u^{\alpha})^t = \overline{u^{\alpha}}^{\ast}$, it
is positive and invertible.  In fact the $i,j$-th element is precisely
$$\sum_{k=1}^{n(\alpha)} \haar ((u^{\alpha}_{ik})^{\ast}u^{\alpha}_{jk})
= n(\alpha) F^{\alpha}_{ij}.$$
It follows that $F^{\alpha}$ is positive, invertible and intertwines
$((u^{\alpha})^t)^{-1}$ and $\overline{u^{\alpha}}$.
\eproof

We will investigate the operators $F^{\alpha}$ further in
section~8.

\bprop
The elements $\{ u^{\alpha}_{pq} \mid \alpha \in I, 1 \leq p,q \leq n(\alpha) \}$
form a basis for $A_0$.
\eprop

\bproof
We still have to show that these elements are linearly independent.
\newline
Let $\sum_{\alpha,i,j}\lambda^{\alpha}_{ij} u^{\alpha}_{ij}$ be any finite
linear combination of the elements $u^{\alpha}_{ij}$, and suppose that
$\sum_{\alpha,i,j}\lambda^{\alpha}_{ij} u^{\alpha}_{ij} =0$.  Let
$\beta \in I$ and $1 \leq r,s \leq n(\beta)$.  Then
$$ 0 = \sum_{\alpha,i,j} \lambda^{\alpha}_{ij}
\haar((u^{\beta}_{rs})^{\ast} u^{\alpha}_{ij})
= \sum_i \lambda^{\beta}_{is} F^{\beta}_{ri}.$$
Because $F^{\beta}$ is invertible, this implies that
$\lambda^{\beta}_{is}=0$ for all $i$. As $s$ and $\beta$ were
arbitrary, this completes the proof.
\eproof

\bprop
We have $\comult (A_0) \subseteq A_0 \odot A_0$ where $A_0 \odot A_0$ is
the algebraic tensor product of $A_0$ with itself.
\eprop

\bproof
Let $v$ be any finite-dimensional unitary representation, acting on a
Hilbert space $\H$.  Fix an orthonormal basis
$\{\xi_1,\, ... \, ,\xi_n \}$ for $\H$.
For any pair of vectors $\xi, \eta \in \H$, we have
\bvgl
\comult ((\omega_{\xi, \eta}\ten \iota)v)
&=& (\omega_{\xi, \eta} \ten \iota \ten \iota)(v_{(12)} v_{(13)})\\
&=& \sum_{i=1}^n (\omega_{\xi, \xi_i} \ten \iota \ten \iota) (v_{(12)})
 (\omega_{\xi_i, \eta} \ten \iota \ten \iota) (v_{(13)})\\
&=& \sum_{i=1}^n (\omega_{\xi, \xi_i} \ten \iota) (v) \ten
 (\omega_{\xi_i, \eta} \ten \iota ) (v).
\evgl
This proves the result.
\eproof

We now want to show that $(A_0, \comult)$ is a Hopf $^{\ast}$-algebra.
We will use the notations 
$\eps$ and $\antipode$ for the
counit and the antipode.  Multiplication in $A_0$ will
sometimes be considered as a linear map $m:A_0 \odot A_0 \rightarrow A_0$
defined by $m(a \ten b) = ab$.  We will also use the summation convention,
which is standard in Hopf algebra theory (see \cite{abe} and \cite{swe})~:
\bvgl
\comult(a) &=& \sum_{(a)} a_{(1)} \ten a_{(2)},\\
\comult^{(2)}(a) &=& (\comult \ten \iota) \comult(a) = \sum_{(a)}
a_{(1)} \ten a_{(2)} \ten a_{(3)},
\evgl
and so on.

In this setting, it is most easy to give the explicit construction
of the counit $\eps$
and the antipode $\antipode$.

Define linear maps $\eps : A_0 \rightarrow  \C$ and
$\antipode: A_0 \rightarrow A_0$ by
\bvgl
\eps (u^{\alpha}_{pq}) &=& \delta_{pq}\\
\antipode(u^{\alpha}_{pq}) &=& (u^{\alpha}_{qp})^{\ast}.
\evgl

\blemma
For all $a \in A$ we have
\bvgl
(\eps \ten \iota)\comult (a) &=& (\iota \ten \eps) \comult(a) = a\\
m(\antipode \ten \iota) \comult(a)
&=& m(\iota \ten \antipode)\comult (a) = \eps (a)1.
\evgl
\elemma

\bproof
For all $\alpha, p,q$, we have
$$(\iota \ten \eps)\comult(u^{\alpha}_{pq}) = \sum_k u^{\alpha}_{pk}
\eps(u^{\alpha}_{kq}) = u^{\alpha}_{pq},$$
so that, by linearity, $(\iota \ten \eps)\comult(a) = a$ for all $a \in A_0$.
Similarly $(\eps \ten \iota)\comult(a) =a$ for all $a$.

Again for all $\alpha,p,q$ we get
$$m(\iota \ten \antipode)\comult(u^{\alpha}_{pq}) = \sum_k u^{\alpha}_{pk}
\antipode(u^{\alpha}_{kq}) = \sum_k u^{\alpha}_{pk}(u^{\alpha}_{qk})^{\ast}
= \delta_{pq} 1 = \eps(u^{\alpha}_{pq})1.$$
And once more by linearity, we have $m(\antipode \ten \iota) \comult(a)
= \eps (a)1$
for all $a \in A_0$.  Similarly $m(\iota \ten \antipode)\comult (a)
= \eps (a)1$
for all $a$.
\eproof

\btheorem
$(A_0, \comult)$ is a Hopf $^{\ast}$-algebra.
\etheorem

\bproof
We have that $A_0$ is a $^{\ast}$-algebra with a comultiplication $\comult$.
Moreover, $\eps:A_0 \rightarrow \C$ is a linear map such that
$$(\eps \ten \iota)\comult = (\iota \ten \eps) \comult = \iota$$
and $\antipode: A_0 \rightarrow A_0$ is a linear map such that
$$m(\antipode \ten \iota) \comult(a) = m(\iota \ten \antipode)\comult (a)
= \eps (a)1.$$
This automatically implies that $\eps$ is a $^{\ast}$-homomorphism, and
that $\antipode$ is an anti-homo\-morphism satisfying
$\antipode(\antipode(a^{\ast})^{\ast})=a$
for all $a \in A$ (see e.g. proposition 2.4
in~\cite{vd.dual}).  So $A_0$ is a Hopf \st-algebra.
\eproof

\bremark \label{dijkhuizen}
In~\cite{dijk}, Dijkhuizen and Koornwinder look at a compact quantum group
through the algebraic object $A_0$.
They define a CQG algebra as a Hopf \st-algebra $A_0$ which is
spanned by the coefficients of its finite-dimensional
unitary (irreducible) corepresentations.
In this setting, the existence of the Haar measure becomes trivial~:
Let $\{u^{\alpha} \mid \alpha \in I \}$ be a complete set
of mutually inequivalent, irreducible unitary representations.
Then the $u^{\alpha}_{pq}$ are linearly independent.  Denote by $1$ the
unique $\alpha$ such that $u^{\alpha}$ is equivalent to the one-dimensional
corepresentation $(1)$.  The linear functional $\haar:A \rightarrow \C$
defined by
$$\haar(u^{\alpha}_{pq}) = \left\{ \begin{array}{ll}
                        0 & \alpha \not= 1 \\
                        1 & \alpha = 1
                        \end{array} \right.$$
is easily shown to be the unique invariant functional on $A_0$
satisfying $\haar(1)=1$.
As a final observation, it is proven that such a CQG
algebra admits a \cst-algebra completion.
However, this completion is not unique.
We will discuss this problem at the end of this section.
\eremark

With the material we have gathered now, it becomes very straightforward
to prove that Woronowicz' new definition of a compact matrix quantum group
in~\cite{wor2} is equivalent with the original one in~\cite{wor1}.

In~\cite{wor1}, a compact matrix quantum group (compact matrix pseudogroup)
is defined to be a pair $(A,u)$ of a unital \cst-algebra $A$ and a matrix
$u$ in $M_n(A)$ satisfying
\begin{enumerate}
\item The \st-subalgebra $A_0$ of $A$, generated by the matrix elements of
$u$, is dense in $A$.
\item There exists a \st-homomorphism $\comult : A \rightarrow A \ten A$
such that
$$\comult(u_{pq})=\sum_k u_{pk} \ten u_{kq}$$
for all $p,q = 1,\ldots,n$.
\item There exists a linear antimultiplicative mapping $\antipode:A_0
\rightarrow A_0$ such that $\antipode(\antipode(a^{\ast})^{\ast})=a$ for all $a \in A_0$
and
\bvgl
\sum_k \antipode(u_{pk})u_{kq} &=& \delta_{pq}1\\
\sum_k u_{pk} \antipode(u_{kq}) &=& \delta_{pq}1
\evgl
for all $p,q = 1,\ldots,n$.
\end{enumerate}
In~\cite{wor2}, this definition is adapted and
a compact matrix quantum group is defined to be a pair $(A,u)$
of a unital \cst-algebra $A$ and a matrix
$u$ in $M_n(A)$ satisfying the axioms 1 and 2 above and the
following condition:
\begin{enumerate}
\item [3'.]
The matrix $(u_{pq})$ and the transpose matrix $(u_{qp})$
are invertible.
\end{enumerate}
Proposition~\ref{prop3.8} says that such a compact matrix
quantum group is a compact quantum group.  Hence, by the previous theorem,
it satisfies also the third of the original axioms above.

Conversely, suppose $(A,u)$ satisfies the original conditions
1, 2 and 3.
By applying $\antipode^{-1}$ to the first equation of axiom 3, we get that
for all $p,q = 1,\ldots,n$
$$\sum_k \antipode^{-1}(u_{kq})u_{pk} = \delta_{pq}1.$$
If we apply the conjugation $^{\ast}$ and use that
$\antipode^{-1}(u_{kq})^{\ast} = \antipode(u_{kp}^{\ast})$, we find
$$\sum_k u_{pk}^{\ast} \antipode(u_{kq}^{\ast}) = \delta_{pq}1.$$
In the same way, by applying $\antipode^{-1}$ to the second equation
of axion 3, we obtain
$$\sum_k \antipode(u_{pk}^{\ast})u_{kq}^{\ast} = \delta_{pq}1.$$
This proves that the matrix $(u_{pq}^{\ast})$ is invertible. Hence
also it's conjugate, the transpose matrix $(u_{qp})$ is invertible.
So $(A,u)$ satisfies axiom 3' as well.
This proves the equivalence
of Woronowicz' original definition of a compact matrix quantum group in
~\cite{wor1}, and his adapted definition in~\cite{wor2}.

We now come back to the faithfulness of the Haar measure.  We saw
in example~\ref{notfaith} that in general, the Haar measure $\haar$ on
$(A, \comult)$ needs not to be faithful. But let us consider the
restriction of $\haar$ to $A_0$.  It gives a Haar measure on the Hopf
algebra $(A_0,\comult)$, which will be faithful.

\bprop
The Haar measure $h$ is faithful on $A_0$.
\eprop

\bproof
Let $a \in A_0$ and suppose that $h(a^\ast a) = 0$.  It follows from the 
Cauchy-Schwarz inequality that then $h(a^\ast b)=0$ for all $b \in A_0$.
We can write $a^\ast =(\rho \ten \iota)v$, where $v$ is a 
finite-dimensional unitary representation of $(A, \Phi)$ on a Hilbert space
$\H$ and $\rho$ is a linear functional on $\B(\H)$.
For any $b \in A_0$ we have 
$$\rho ((\iota \ten h)(v(1 \ten b)) = h(a^\ast b) = 0.$$
If we denote the \cst-subalgebra of $\B(\H)$ defined in
proposition~\ref{propos} by $\B$, then we obtain $\rho (\B) = 0$.
Because $v \in \B \ten A$
this gives us $a^\ast = (\rho \ten \iota)v = 0$.
\eproof

In fact we have the following more general property (proposition~3.4
from~\cite{vd.alg}): Let $(A_0, \comult)$ be any Hopf algebra and $\haar$
a non-zero linear functional on $A_0$ which is left invariant.
Then $\haar$ is left faithful
in the sense that $\haar(ba) = 0$ for all $b \in A_0$ implies $a=0$.
If $(A_0, \comult)$ is a Hopf algebra with an invertible antipode, we also
have faithfulness on the other side.  
In the case we consider, $(A_0, \comult)$ is a Hopf \st-algebra,
so $\antipode$ is automatically invertible.
If moreover $\haar$ is positive, then left or right faithfulness
in the above sense are
of course equivalent with the usual notion~: $\haar(a^{\ast}a)=0 \Rightarrow
a=0$.

On the \cst-algebra level, we have the following property~:

\bprop
Let $a \in A$.  Then $\haar(a^{\ast}a)=0$ if and only if $\haar(aa^{\ast})=0$.
\eprop

\bproof
We will use the result of proposition 3.12 of \cite{vd.alg}.  It says that
if $A_0$ is a regular multiplier Hopf algebra with a left invariant functional
$\haar$, then there is an automorphism $\sigma$ of $A_0$ such that $\haar(ab)=
\haar(b \sigma(a))$ for all $a,b \in A_0$.  In our case, $A_0$ is a Hopf
\st-algebra and hence automatically regular.

We only prove one implication, the other one is completely analogous.
Suppose that $a \in A$ and that $\haar(a^{\ast}a)=0$.  Then by
the Cauchy-Schwarz inequality, $\haar(ba)=0$ for all $b \in A$.
Using the result mentioned above and the continuity of $\haar$, we obtain
$\haar(ba)= \haar(a \sigma(b)) =0$
for all $b \in A_0$.
Now we can approximate $a^{\ast}$ by linear combinations of elements of
the form $\sigma(b)$ with $b \in A_0$.  So we obtain that
$\haar(aa^{\ast})=0$.
\eproof

So $I=\{a \in A \mid \haar(a^{\ast}a) =0 \}$ is a closed two-sided
ideal in $A$.
Hence we can consider the quotient $A_r=A/I$.  It is clear that
$(A_r, \comult)$ is again a compact quantum group, where now the Haar
measure is faithful.  Of course, the associated Hopf \st-algebra is still
the same.  So the process of taking the quotient amounts to considering
another \cst-norm on $A_0$ and taking the completion with respect to this
norm.  This also means that both these compact quantum groups have the same
representations, which will imply that they have the same dual (see
section~8).  These observations suggest that we should consider
$(A, \comult)$ and $(A_r, \comult)$ as the same object, as different forms
of the same compact quantum group.

\section{The dual of a compact quantum group}

Let $G$ be an abelian locally compact group.  The set of characters of $G$
can be made into an abelian locally compact group in a natural way.  It is
called the dual group of $G$, and the Pontryagin-Van Kampen duality theorem
says that its dual is again (isomorphic to) the original group $G$.  For
compact groups, the Tannaka-Krein duality generalises this duality theory
to the non-abelian case.  Several authors have discussed a duality theory
in the quantum framework: see for instance 
\cite{baaj, es, hol, hol2, mn, m, mnw, 
vd.alg, vd.alg2, wang.krein, wor.tan}.
In these theories, the unitary irreducible representations play the role
the characters have in the case of an abelian group.  The idea is that the
representations contain sufficient information to fully recover the (quantum)
group.

In this section, we will show how the irreducible unitary representations of
a compact quantum group give rise to a discrete quantum group in the sense
of~\cite{vd.discr}.  We want to mention that this duality between compact
and discrete quantum groups is a special case of the duality theory
developed in~\cite{vd.alg} and~\cite{vd.alg2}.

Let $(A, \comult)$ be a compact quantum group.
To construct the dual, we start
with a complete set $\{u^{\alpha} | \alpha \in I \}$ of mutually
non-equivalent, irreducible unitary representations.
Let $A_0$ be the associated Hopf \st-algebra as in section~7.
Define $B_0$ to be the space of
linear functionals on $A$ defined by $x \mapsto \haar(ax)$ where $a \in A_0$.
So $B_0$ is a subspace of the dual $A^{\ast}$.
For $\alpha \in I$, let $F^{\alpha}$ be the positive invertible operator
obtained in proposition~\ref{orth}.

\blemma
Let $$a= \sum_{k=1}^{n(\alpha)} ((F^{\alpha})^{-1})_{pk}
(u^{\alpha}_{kq})^{\ast}$$
and let $\omega^{\alpha}_{pq}$ be the linear functional on $A$ defined
by $\omega^{\alpha}_{pq}(x) = \haar(ax)$.  Then
$$\omega^{\alpha}_{pq} (u^{\alpha}_{pq})=1$$
and $\omega^{\alpha}_{pq}$ is $0$ on all other matrix elements.
\elemma

\bproof
We have
$$\omega^{\alpha}_{pq} (u^{\beta}_{rs}) = \sum_{k=1}^{n(\alpha)}
((F^{\alpha})^{-1})_{pk} \haar((u^{\alpha}_{kq})^{\ast} u^{\beta}_{rs})
= \delta_{\alpha \beta} \delta_{qs} \sum_{k=1}^{n(\alpha)}
((F^{\alpha})^{-1})_{pk} F^{\alpha}_{kr}\\
= \delta_{\alpha \beta} \delta_{qs} \delta_{pr}.$$

\eproof

It is clear that $B_0$ is spanned by the elements $\{ \omega^{\alpha}_{pq}
\mid \alpha \in I, 1 \leq p,q \leq n(\alpha) \}$.

As before, we consider the product on $A^{\ast}$ induced by the
comultiplication on $A$.  This is defined by $(\omega \psi)(a)
= (\omega \ten \psi) \comult (a)$.  We have
\bvgl
\omega_{pq}^{\alpha} \omega_{rs}^{\beta} (u^{\gamma}_{ij}) &=&
\sum_l \omega_{pq}^{\alpha} (u_{il}^{\gamma}) \omega_{rs}^{\beta}
(u_{lj}^{\gamma})\\
&=& \delta_{\alpha \gamma} \delta_{\beta \gamma} \delta_{pi}
\delta_{sj} \delta_{qr} = \delta_{\alpha \beta} \delta_{qr}
\omega^{\alpha}_{ps}(u^{\gamma}_{ij}).
\evgl
It follows that $B_0$ is a subalgebra of $A^{\ast}$ and that this
subalgebra is the algebraic direct sum $\sum_{\alpha \in I} \oplus
M_{n(\alpha)}$.

The antipode $\antipode$ induces a \st-structure on $B_0$ by
$$(\omega^{\alpha}_{pq})^{\ast}(x) = \omega^{\alpha}_{pq}
(\antipode(x)^{\ast})^-$$
for $x \in A_0$.
Now, as $\antipode(u^{\alpha}_{pq}) = (u^{\alpha}_{qp})^{\ast}$,
we get $(\omega^{\alpha}_{pq})^{\ast} = \omega^{\alpha}_{qp}$.
So we get the
usual \st-algebra structure on $B_0$.
Remark that $B_0$ has no identity:  $\eps(u^{\alpha}_{pq}) =
\delta_{pq}$ but $\eps$ does not belong to $B_0$.

Now we want to make $B_0$ into a discrete quantum group as
in~\cite{vd.discr}.  Let us recall the definition:

\bdef
A discrete quantum group is a pair $(B_0, \comult)$ where $B_0$
is a direct sum
of full matrix algebras and $\comult$ is a comultiplication on $B_0$ making
$B_0$ into a multiplier Hopf \st-algebra.
\end{definition}

So we have to make $B_0$ into a multiplier Hopf  \st-algebra (see
definitions 2.3 and 2.4 in~\cite{vd.mult}).
The multiplier algebra of $B_0 \odot B_0$ consists of course of
all elements in $\prod_{\alpha, \beta} M_{n(\alpha)}\ten M_{n(\beta)}$.
Define $\comult$ on $B_0$ by
$$\comult (\omega)(a \ten b) = \omega (ab),$$
then we get a linear map from $B_0$ into this multiplier algebra
$M(B_0 \odot B_0)$.  It is automatically a \st-homomorphism~:
\bvgl
\comult (\omega \psi)(a \ten b) &=& (\omega \psi)(ab)\\
&=& (\omega \ten \psi) \sum_{(ab)} (ab)_{(1)} \ten (ab)_{(2)}\\
&=& \sum_{(a)(b)} \omega(a_{(1)}b_{(1)}) \psi(a_{(2)}b_{(2)})\\
&=& \sum_{(a)(b)} \comult(\omega)(a_{(1)} \ten b_{(1)})
\comult(\psi)(a_{(2)} \ten b_{(2)})\\
&=& \comult(\omega) \comult(\psi)(a \ten b).
\evgl
Similarly,
\bvgl
\comult (\omega^{\ast})(a \ten b) &=& \omega^{\ast}(ab)\\
&=& \omega (\antipode(ab)^{\ast})^- \\
&=& \omega (\antipode(a)^{\ast}\antipode(b)^{\ast})^- \\
&=& \comult(\omega) (\antipode(a)^{\ast} \ten \antipode(b)^{\ast})^-\\
&=& \comult(\omega)^{\ast} (a \ten b).
\evgl

Now we can define the antipode.  We know that
$\antipode(u^{\alpha}_{pq})
= ({u^{\alpha}_{qp}})^{\ast}$ and that for $\alpha \in I$
there is $\overline{\alpha} \in I$ such that
$\overline{u^{\alpha}}$ is
unitarily equivalent with $u^{\overline{\alpha}}$.  Hence $\antipode$ maps
$B_0$ to $B_0$.
In fact $\antipode(B_{\alpha}) \subseteq B_{\overline{\alpha}}$
and $\antipode(B_{\overline{\alpha}}) \subseteq B_{\alpha}$.
So we could use the antipode $\antipode$ on $A$
to define the antipode on the dual, and prove in this way that $B_0$ is
a multiplier Hopf algebra.
But the following method is easier.

Recall that, because $(A_0, \comult)$ is a Hopf algebra, we have that
\bvgl
T_1 : x \ten y &\mapsto& \comult(x)(1 \ten y)\\
T_2 : x \ten y &\mapsto& (x \ten 1)\comult(y)
\evgl
extend to bijections from the algebraic tensor product of $A_0 \odot A_0$
to itself.  Now, for $\omega, \psi$ in $B_0$ we have
\bvgl
((\omega \ten 1)\comult(\psi))(x \ten y) &=& \sum_{(x)} \omega (x_{(1)})
\comult(\psi)(x_{(2)} \ten y)\\
&=& \sum_{(x)} \omega (x_{(1)})\psi(x_{(2)}  y)\\
&=& (\omega \ten \psi)(\comult(x)(1 \ten y)).
\evgl
Let $a,b \in A_0$ such that $\omega$ and $\psi$ are the functionals
given by $x \mapsto \haar(ax)$ and $x \mapsto \haar(bx)$ respectively.
As we have already mentioned in section~7, the Hopf algebra $A_0$ is
automatically regular, because it is a Hopf \st-algebra.  This means
that the antipode is invertible, which is equivalent to the condition
that $(A_0,\comult')$ is also a Hopf algebra, where $\comult'$ denotes
the opposite comultiplication.
So we can write $a \ten b = \sum (1 \ten p_i)\comult(q_i)$
with $p_i, q_i \in A_0$.  Now
\bvgl
((\omega \ten 1)\comult(\psi))(x \ten y)
&=& (\haar \ten \haar)((a \ten b)\comult(x)(1 \ten y))\\
&=& (\haar \ten \haar)(\sum (1 \ten p_i)\comult(q_i x)(1 \ten y))\\
&=& \sum \haar(q_ix) \haar(p_iy).
\evgl
It follows that $(\omega \ten 1)\comult(\psi) \in B_0 \odot B_0$.
Similarly,
$$(\comult(\omega)(1 \ten \psi))(x \ten y) = (\omega \ten \psi)((x \ten 1)
\comult(y))$$
and $\comult(\omega)(1 \ten \psi) \in B_0 \odot B_0$. Here we need also
the fact that for a regular (multiplier) Hopf algebra $A_0$ with left
invariant functional $h$, the sets
$\{h(a\,\cdot\,) \mid a \in A_0 \}$ and $\{h(\,\cdot\,a) \mid a \in A_0 \}$
are equal (proposition 3.11 in \cite{vd.alg}).

In the same way, the inverses of the maps $T_1, T_2$ have duals; and these
are the inverses of the duals of $T_1$ and $T_2$.
This implies that $(B_0, \comult)$ is a multiplier Hopf algebra.

As $B_0$ is a direct sum of full matrix algebras,
the \cst-norm on $B_0$ exists and is unique.
The completion $B$ gives the discrete quantum group
on the \cst-algebra level.

The element $h$ in proposition 3.1 of~\cite{vd.discr} is of course the Haar
functional $\haar$ itself.  Let $\tr_{\alpha}$ denote the trace on $B_{\alpha}$
and consider it as a linear map on all of $B_0$.  It is proven
in~\cite{vd.discr} that the
operators $K_{\alpha} = (\tr_{\alpha} \ten \iota)\comult(h)$ implement
$\antipode$,
i.e.\  $\antipode^2(\omega) = K_{\alpha}^{-1} \omega K_{\alpha}$
for all $\omega \in B_{\alpha}$.
What are these operators here~? We have
that $\tr_{\alpha} = \tr_{\overline{\alpha}} \circ \antipode$.  So
\bvgl
\langle (\tr_{\alpha} \ten \iota)\comult(h), u^{\beta}_{pq} \rangle
&=& \langle (\tr_{\overline{\alpha}} \circ \antipode \ten \iota)\comult(h),
u^{\beta}_{pq} \rangle\\
&=& \sum_{k=1}^{n(\alpha)} \haar((u^{\overline{\alpha}}_{kk})^{\ast}
u^{\beta}_{pq})
= \delta_{\beta \overline{\alpha}} F^{\overline{\alpha}}_{qp}.\\
\evgl

In this way one can construct the dual quantum group out of the
unitary representations of a compact quantum group.
Conversely, starting from a discrete quantum group, one can recover the Hopf
\st-algebra $A_0$.  The \cst-norm on $A_0$ is not unique, but as we
discussed at the end of section~7, we could consider the different possible
\cst-completions as the same compact quantum group.

It is clear that this duality generalizes the Pontryagin duality between
compact and discrete abelian groups: If $A=C(G)$ is the function algebra
of a compact abelian group, then the irreducible unitary representations
are given by the characters of $G$.  The dual discrete quantum group $B_0$
is then the function algebra of the character group.

\end{document}